\documentclass[a4paper,12pt,reqno]{amsart}
\usepackage{amssymb}
\usepackage{amsmath}
\usepackage{hyperref}
\usepackage{mathrsfs}
\usepackage{bbm}
\usepackage{tikz}
\usetikzlibrary{commutative-diagrams}
\usepackage{tikz-cd}
\usepackage{bm}
\usepackage{bm}
\usepackage{geometry}
\usepackage{colordvi}
\usepackage{color}

\allowdisplaybreaks

\numberwithin{equation}{section}

\theoremstyle{definition}
\newtheorem{theorem}{Theorem}[section]
\newtheorem{corollary}[theorem]{Corollary}
\newtheorem{lemma}[theorem]{Lemma}
\newtheorem{definition}[theorem]{Definition}
\newtheorem{proposition}[theorem]{Proposition}
\newtheorem{example}[theorem]{Example}
\newtheorem{remark}[theorem]{Remark}

\newtheorem*{THM}{Theorem}

\begin{document}

\title{Twisted Poincar\'e duality for orientable Poisson manifolds}

\author{Tiancheng Qi}
\address{School of Mathematical Sciences, Fudan University, Shanghai 200433, China}
\email{tcqi21@m.fudan.edu.cn}

\author{Quanshui Wu}
\address{School of Mathematical Sciences, Fudan University, Shanghai 200433, China}
\email{qswu@fudan.edu.cn}



\begin{abstract}
We geometrize the constructions of twisted Poisson modules introduced by Luo-Wang-Wu, and Poisson chain complexes with coefficients in Poisson modules defined in the algebraic setting to the geometric setting of Poisson manifolds. We then prove that for any orientable Poisson manifold $M$, there is an explicit chain isomorphism between the Poisson cochain complex with coefficients in any Poisson geometric module and the Poisson chain complex with coefficients in the corresponding twisted Poisson geometric module, induced by a modular vector field of $M$. These are the geometric analogues of results obtained by Luo-Wang-Wu for smooth Poisson algebras with trivial canonical bundle. In particular, a version of twisted Poincar\'e duality is established between the Poisson homologies and the Poisson cohomologies of an orientable Poisson manifold with coefficients in an arbitrary vector bundle with a flat contravariant connection. This generalizes the duality theorems for orientable Poisson manifolds established by Evens-Lu-Weinstein, and by Xu.
\end{abstract}


\thanks{This research has been supported by the NSFC (Grant No. 12471032) and the National Key Research and Development Program of China (Grant No. 2020YFA0713200).}

\maketitle	

\section*{Introduction}
Poisson structures, rooted in mechanics, have nourished many branches of mathematics and physics over the past decades, such as classical/quantum mechanics, non-commutative geometry/algebra, and representation theory. Poisson structures on smooth manifolds were introduced by Lichnerowicz \cite{MR0501133}. Any symplectic manifold admits a canonical Poisson structure, but a Poisson structure on a manifold does not necessarily come from a symplectic structure. 

The (co)homology theory of Poisson structures is crucial in the study of Poisson structures and their deformations. Poisson cohomology spaces carry important information about Poisson structures. For example, the space of all Casimir functions on a Poisson manifold is the 0th cohomology; Poisson vector fields modulo Hamiltonian vector fields is the 1st cohomology; and the 2nd and 3rd cohomologies control the formal deformations of the Poisson structure. For any orientable Poisson manifold, there exists a special characteristic class in its first cohomology, known as the modular class, which was introduced by Weinstein \cite{MR1484598}. This cohomology class measures the existence of a volume form that is preserved by all Hamiltonian flows. Poisson cohomology also plays a fundamental role in the study of the normal forms of Poisson structures, see \cite{MR2178041} for a systematic discussion. Poisson cohomology of a symplectic manifold coincides with its de Rham cohomology \cite{MR0501133}. In general, the computation of both homology and cohomology for a given Poisson structure is not easy. There are deep relations between Poisson (co)homology and Hochschild (co)homology; see, for example, \cite{MR950556,MR922799,MR2062626}. The Hochschild homology and cyclic homology of some non-commutative algebras can be calculated by  viewing these algebras as deformations of Poisson algebras and employing the Brylinski spectral sequence \cite{MR950556}. Conversely, the Poisson homology of a Poisson structure can also be computed by using the Hochschild homology of its deformation algebra, see, for example, \cite{MR1291019,MR2485483,MR2301394,MR3254768}. Dolgushev proved that for any smooth Poisson affine variety with trivial canonical bundle, the Poisson structure is unimodular if and only if its (formal) deformation quantization algebra is Calabi-Yau \cite{MR2480714}. L\"u-Wang-Zhuang showed that the Poisson enveloping algebra of any smooth affine Poisson algebra with trivial canonical bundle is skew Calabi-Yau, and when the Poisson structure is unimodular, the enveloping algebra is Calabi-Yau \cite{MR3687261}. For a more general version of this result in the context of the universal enveloping algebra of Lie-Rinehart algebras, see \cite{MR3893434}. Wu-Zhu considered filtered deformations of Poisson structures and proved that the deformation algebra is Calabi-Yau if and only if the corresponding Poisson structure is unimodular under some mild assumptions \cite{MR4200729}. In fact, the relation was described in \cite{MR4200729} between the Nakayama automorphisms of the filtered deformations and the modular derivations of the Poisson structures by using homological determinants as a bridge.

On the geometric side, the notion of Poisson cochain complexes and Poisson cohomologies of a Poisson manifold were introduced by Lichnerowicz \cite{MR0501133}. The definitions of  Poisson chain complexes and Poisson homologies of a Poisson manifold, which can be traced back to the work of Gel'fand-Dorfman \cite{MR684122}, Koszul  \cite{MR0837203} and Brylinski \cite{MR950556}, are also referred to as  canonical complexes and canonical homologies by Brylinski, respectively. Poisson cohomology is also studied in the more general setting of Lie algebroids, introduced by Pradines \cite{MR216409}. For any Poisson manifold, its cotangent bundle has a natural Lie algebroid structure \cite {MR996653,MR959095}. Hence, one can consider Poisson cohomology with coefficients as in the Lie algebroid case \cite{MR896907}. In \cite{MR950556}, a version of Poincar\'e duality was  established  between the Poisson homology and the Poisson cohomology of symplectic manifolds. In general, there exist Poisson manifolds for which there is no Poincar\'e duality between the Poisson homology and the Poisson cohomology, see \cite[Example 5.19]{MR1269545} for example. For any Poisson manifold that admits a nowhere vanishing density, Evens-Lu-Weinstein constructed an explicit representation on the canonical line bundle and established a ``twisted Poincar\'e duality'' between Poisson homology and the Poisson cohomology with coefficients in the canonical line bundle \cite{MR1726784}. Xu also discovered the representation on the canonical bundle of Poisson manifolds constructed by Evens-Lu-Weinstein, and he proved that there is a Poincar\'e duality between Poisson homology and Poisson cohomology for any orientable unimodular Poisson manifold \cite{MR1675117}.  Sti\'enon \cite{MR2764872} established holomorphic analogues of the aforementioned results by Evens-Lu-Weinstein and Xu for holomorphic Poisson manifolds. For a discussion with a geometric point of view on Poincar\'e duality in the context of Lie algebroids, see \cite{MR1726784,MR1675117,MR2369386,MR2228680} for a partial list.

On the algebraic side, Huebschmann constructed a natural Lie-Rinehart algebra structure on the module of K\"ahler differentials of a Poisson algebra, which can be regarded as an algebraic incarnation of the Lie algebroid structure on the cotangent bundle of a Poisson manifold \cite{MR1058984}. With the aid of the homological tools developed in \cite{MR0154906} by Rinehart, Huebschmann defined the Poisson cohomology and homology with coefficients for Poisson algebras and also defined geometric Poisson cohomology and homology with coefficients for Poisson manifolds in an algebraic way. He also proved that when restricted to the context of Poisson manifolds, algebraic Poisson cohomology coincides with geometric Poisson cohomology. However, this is not the case for Poisson homology. In \cite{MR1696093}, a general duality theorem is proved by Huebschmann for any duality Lie-Rinehart algebra with constant rank, which covered the twisted Poincar\'e duality given by Evens-Lu-Weinstein. For any Poisson module over a Poisson algebra (to our knowledge, the first paper to provide a formal definition of Poisson module were \cite{MR1423640} by Reshetikhin-Voronov-Weinstein, although they referred to it as Poisson bimodule), Fresse \cite{MR1647186,MR2275449} explicitly formulated the (algebraic) Poisson chain complex with coefficients in any Poisson module (see also \cite{maszczyk2006}), which is consistent with Huebschmann's definition in the setting of Lie-Rinehart algebras. Fresse also used two different algebraic approaches to define Poisson homology and proved that in the case of smooth Poisson algebras, the homology defined by these two methods coincides with the homology given by the Poisson chain complex. Knowing that Poisson (co)homology has a close relation with the Hochschild (co)homology of its deformation algebra, Launois-Richard \cite{MR2301394} obtained a new Poisson module by considering the semiclassical limit of the dualizing bimodule between Hochschild homology and cohomology of the quantum affine space, which established a twisted Poincar\'e duality for some quadratic Poisson algebras. Using similar methods, Zhu achieved a twisted Poincar\'e duality for some linear Poisson algebras \cite{MR3314106}. In \cite{MR3395070}, Luo-Wang-Wu generalized the aforementioned results to any polynomial Poisson algebra. The key factor in achieving the duality is that any Poisson module (over a polynomial Poisson algebra) can be naturally twisted by a modular derivation, which can be regarded as the algebraic version of the modular vector field in Poisson geometry. Later, L\"u-Wang-Zhuang provided a twisted Poincar\'e duality for any affine smooth Poisson algebra by using the Serre invertible bimodule \cite{MR3687261}. Luo-Wang-Wu extended the ideas from \cite{MR3395070}, and achieved a twisted Poincar\'e duality at the chain level for any affine smooth Poisson algebra with trivial canonical bundle \cite{MR4727073}. Chemla in \cite{MR4761654} proved a twisted Poincar\'e duality for Hopf algebroids with bijective antipode, and thereby recovered Huebschmann's duality theorem particularly for Poisson manifolds. See also \cite{MR3264266,MR3969113,MR4127112,MR4345205,MR4575095} and references therein for related works on twisted Poincar\'e duality for Poisson algebras.

Note that Poisson chain complexes with coefficients, and twisted Poincar\'e duality with arbitrary coefficients between Poisson homology and cohomology are discussed only in the algebraic realm. It will be interesting to geometrize some of the constructions and ideas in Poisson algebras.

To warm up, let us consider the geometrization of  twisted Poisson modules (by a Poisson derivation) introduced by Luo-Wang-Wu \cite{MR3395070}. In the rest of the introduction, let $(M,\pi)$ be a Poisson manifold of dimension $n$, $E$ be a smooth vector bundle over $M$ with a flat $T^{*}M$-connection on $E$, say $\nabla$. Denote the algebra of smooth functions on $M$ by $C^{\infty}(M)$, and the space of all smooth global sections of $E$ by $\Gamma(E)$, respectively. Then for any Poisson vector field $\phi$ on $M$, we can twist $\nabla$ by the vector field $\phi $ to obtain a new flat $T^{*}M$-connection $\nabla^{\phi}$ on $E$ (see Proposition \ref{prop-twist-flatc}) such that
$$\nabla_{df}^{\phi }(s)=\nabla_{df}(s)-\phi(f)s,\forall f\in C^{\infty}(M),s\in\Gamma(E).$$
When $M$ is orientable with a volume form $\mu$, let $\phi_{\mu}$ be the modular vector field with respect to $\mu$, and let $\nabla_{0}$ be the obvious flat $T^{*}M$-connection on the canonical line bundle $\wedge^{n}T^{*}M$ given by $(\nabla_{0})_{\omega}\mu=0,\forall \omega\in\Omega^{1}(M)$, then the twisted flat $T^{*}M$-connection $(\nabla_{0})^{\phi_{\mu}}$ of $\nabla_{0}$ by $\phi_{\mu}$ is given by the formula
$$(\nabla_{0})_{\omega}^{\phi_{\mu}}s=\omega\wedge d\iota_{\pi}(s),\forall \omega\in\Omega^{1}(M),s\in \Omega^{n}(M),$$
where $\iota_{\pi}$ denotes the interior product by $\pi$. This is the canonical flat $T^{*}M$-connection on $\wedge^{n}T^{*}M$ discovered by Evens-Lu-Weinstein \cite{MR1726784} and Xu \cite{MR1675117} respectively.

In this paper, we geometrize Poisson chain complexes with coefficients in an arbitrary Poisson module in the context of Poisson algebra. The critical ingredient of this geometrization is the terminology of the geometric module (see Definition \ref{geo-module-def}), introduced by Krasil'shchik-Lychagin-Vinogradov \cite{MR0861121}. A geometric $C^{\infty}(M)$-module equipped with a Poisson module structure (over the Poisson algebra $C^{\infty}(M) $) is referred to as a Poisson geometric module, which is, in fact, a generalization of a vector bundle with a flat $T^{*}M$-connection (see Proposition \ref{geomodule-poisson-corresprop} and Remark \ref{rmk-poimodule-notE}). Suppose $(W,\{-,-\}_{W})$ is a right Poisson geometric module over $C=C^{\infty}(M)$, then the Poisson chain complex of $M$ with coefficients in $W$ (see Proposition \ref{prop-poisson-chain}) can be formulated as
\begin{equation}\label{intro-poiss-complex}
\begin{tikzcd}[column sep=small]
\cdots \arrow[r] & W \otimes_{C} \Omega^q(M) \arrow[r, "\partial_q"] & W \otimes_C \Omega^{q-1}(M) \arrow[r] & \cdots \arrow[r, "\partial_2"] & W \otimes_C \Omega^1(M) \arrow[r, "\partial_1"] & W \arrow[r] & 0,
\end{tikzcd}
\end{equation}
where $\partial_{r}: W\otimes_{C}\Omega^{r}(M)\to W\otimes_{C}\Omega^{r-1}(M)$ is defined as:
\begin{align*}
\partial_{r}(w\otimes f_{0}df_{1}\wedge\cdots\wedge df_{r})=&\sum\limits_{i=1}^{r}(-1)^{i-1}\{wf_{0},f_{i}\}_{W}\otimes df_{1}\wedge\cdots\widehat{df_{i}}\cdots\wedge df_{r}\\
&+\sum\limits_{1\leq i < j\leq r}(-1)^{i+j}w\otimes f_{0}d\{f_{i},f_{j}\}\wedge df_{1}\wedge\cdots\widehat{df_{i}}\cdots\widehat{df_{j}}\cdots\wedge df_{r}.
\end{align*}

The difference between complex (\ref{intro-poiss-complex}) and the algebraic Poisson chain complex (see (\ref{alg-Poisson-complex})) is that the modules of differentials in (\ref{intro-poiss-complex}) consist of smooth forms, while those in the definition of the algebraic Poisson chain complex consist of K\"ahler forms.


Since a flat $T^{*}M$-connection on $E$ corresponds to a natural Poisson module structure on the geometric module $\Gamma(E)$ \cite{MR1959580,MR2021638},
one also obtains the Poisson chain complex of $M$ with coefficients in $(E,\nabla)$ (see Corollary \ref{cor-poissoncomplex-flat}). In this case, the homology of the complex (\ref{intro-poiss-complex}) coincides with the definition of the geometric Poisson homology with coefficients in $(E,\nabla)$ given by Huebschmann (see Remark \ref{rmk-coinc-Hueb-ho}). When $W=C^{\infty}(M)$ with the obvious Poisson module structure, the complex (\ref{intro-poiss-complex}) boils down to the canonical complex of $M$. Since $\Omega^{1}(M)$ is a projective $C^{\infty}(M)$-module, the complex (\ref{intro-poiss-complex}) can also be derived from a more general construction by Rinehart in the context of Lie-Rinehart algebras \cite{MR0154906}. Our construction here is more direct and simpler, obtained by appropriately adjusting the differentials of the canonical complex of $M$.

The major work in this paper is to prove the geometric version of the twisted Poincaré duality for smooth Poisson algebras with trivial canonical bundle established by Luo-Wang-Wu \cite{MR4727073}. In algebraic  geometry, the triviality of the canonical bundle of an affine smooth variety means that the variety is Calabi-Yau, while in differential geometry, the triviality of  the canonical bundle of a smooth manifold means that the manifold is orientable. This motivates the following result.
\begin{THM}
Let $M$ be an orientable Poisson manifold of dimension $n$ with a volume form $\mu$, and let $\phi_{\mu}$ be the modular vector field of $M$ with respect to $\mu$.
\begin{itemize}
\item[(1)]For any Poisson geometric module $W$ and its twisted Poisson module $W_{-\phi_{\mu}}$, there is an explicit isomorphism between the Poisson cochain complex of $M$ with coefficients in $W$ and the Poisson chain complex of $M$ with coefficients in $W_{-\phi_{\mu}}$. In particular, $\text{HP}^{k}(M,W ) \cong \text{HP}_{n-k}(M,W_{-\phi_{\mu}}),\  0\leq k\leq n,$
where $\text{HP}^{k}(M,W)$ is the $k$-th Poisson cohomology of $M$ with coefficients in $W$, and $\text{HP}_{n-k}(M,W_{-\phi_{\mu}})$ is the $(n-k)$-th Poisson homology of $M$ with coefficients in the twisted Poisson geometric module $W_{-\phi_{\mu}}$. Moreover, $\text{HP}_{n-k}(M,W)\cong \text{HP}^{k}(M,W_{\phi_{\mu}}),\  0\leq k\leq n$.
\item[(2)]For any vector bundle $E$ over $M$ with a flat $T^{*}M$-connection $\nabla$ and its twisted flat $T^{*}M$-connection $\nabla^{-\phi_{\mu}}$, there is an explicit isomorphism between the Poisson cochain complex of $M$ with coefficients in $(E,\nabla)$ and the Poisson chain complex of $M$ with coefficients in $(E,\nabla^{-\phi_{\mu}})$. In particular, $\text{HP}^{k}(M,\nabla) \cong \text{HP}_{n-k}(M,\nabla^{-\phi_{\mu}}),\  0\leq k\leq n,$
where $\text{HP}^{k}(M,\nabla)$ is the $k$-th Poisson cohomology of $M$ with coefficients in $(E,\nabla)$, and $\text{HP}_{n-k}(M,\nabla^{-\phi_{\mu}})$ is the $(n-k)$-th Poisson homology of $M$ with coefficients in $(E,\nabla^{-\phi_{\mu}})$. Moreover, $\text{HP}_{n-k}(M,\nabla)\cong \text{HP}^{k}(M,\nabla^{\phi_{\mu}}),\  0\leq k\leq n$.
\end{itemize}
\end{THM}
Since the aforementioned canonical flat $T^{*}M$-connection on $\wedge^{n}T^{*}M$ can be realized by twisting the obvious flat $T^{*}M$-connection induced by the given volume form, we recover the twisted Poincar\'e duality theorem established by Evens-Lu-Weinstein for orientable Poisson manifolds \cite{MR1726784}. Our result also implies that there is a Poincar\'e duality between Poisson homologies and Poisson cohomologies for any orientable unimodular Poisson manifold $M$ with coefficients in an arbitrary flat $T^{*}M$-connection, which generalizes the duality theorem of Xu \cite{MR1675117}. 
\par There are general duality theorems connecting the Poisson homologies and the Poisson cohomologies of Poisson manifolds which were proved by Huebschmann \cite{MR1696093} in the context of Lie-Rinehart algebras, and by Chemla \cite{MR4761654} in the setting of Hopf algebroids, respectively. Our approach contrasts with Huebschmann's in the following respect. Retain both the assumptions and notations of the above theorem. For any Poisson geometric module $W$, the construction of the twisted Poisson module introduced in \cite{MR3395070} allows us to modify a family of standard isomomorphisms from each term of the Poisson cochain complex with coefficients in $W$ to the corresponding term of the Poisson chain complex with coefficients in $W$, defined by the interior product and a fixed volume form, into a chain isomorphism from the Poisson cochain complex with coefficients in $W$ to the Poisson chain complex with ``twisted coefficients'' in $W_{-\phi_{\mu}}$. This relies on some useful formulas given by Luo-Wang-Wu in the context of Poisson algebras \cite{MR3395070,MR3969113,MR4727073}. Our treatment retains the Poisson cochain complex with coefficients, which in special cases is the Chevalley-Eilenberg cochain complex for computing the cohomology of Poisson Lie algebroids with coefficients in a representation. Huebschmann established a duality theory 
by employing derived functors and the universal enveloping algebra of a Lie-Rinehart algebra. Compared with Chemla's result, the isomorphisms between Poisson homologies and Poisson cohomologies established in this paper are induced by an explicit chain isomorphism between the Poisson cochain complex with coefficients and the Poisson chain complex with twisted coefficients, whereas Chemla proved that there is an isomorphisms between the homologies and the cohomologies by working in the derived category.

This paper is organized as follows. In Section \ref{section 1}, we collect some facts concerning K\"ahler differentials, multi-derivations, multivector fields, and geometric modules in the sense of Krasil'shchik-Lychagin-Vinogradov \cite{MR0861121}. In Section \ref{section 2}, after reviewing some basic concepts about Lie algebroids, Lie-Rinehart algebras, and Poisson algebras, we proceed to establish a correspondence between flat $\Omega^{1}(M)$-connections, flat $\Omega_{\mathbb{R}}^{1}(C^{\infty}(M))$-connections, right Poisson module structures, and right $\mathcal{U}(C^{\infty}(M),\Omega^{1}(M))$-module structures that are compatible with the fixed $C^{\infty}(M)$-module structure on any geometric $C^{\infty}(M)$-module (see Propositions \ref{geomodule-poisson-corresprop} and \ref{prop-RLmodule-right}). Here, $\mathcal{U}(C^{\infty}(M),\Omega^{1}(M))$ denotes the universal enveloping algebra of the Lie-Rinehart algebra $(C^{\infty}(M),\Omega^{1}(M))$. In Section \ref{section 3}, we geometrize the algebraic Poisson chain complex with coefficients in any Poisson module and the twisted Poisson module introduced in \cite{MR3395070}, respectively, to the geometric setting of Poisson manifolds (see Propositions \ref{prop-poisson-chain} and \ref{prop-twist-flatc}). Subsequently, we establish the twisted Poincar\'e duality between the Poisson homologies and the Poisson cohomologies with coefficients for orientable Poisson manifolds (see Theorems \ref{thm-poincare-dual-chain}, \ref{thm-poincare-dual-homo} and \ref{thm-poincare-flatconnver}).

\section{Preliminaries}\label{section 1}

In this section, we collect some basic facts regarding K\"ahler differentials, multi-derivations, geometric modules, multivector fields, and interior products. Throughout this paper, all manifolds considered are real, smooth, Hausdorff, and second countable. Additionally, all morphisms between manifolds, as well as vector bundles over these manifolds, are assumed to be smooth.  Let $M$ be a smooth manifold. If $E$ is a vector bundle over $M$, the space of all smooth global sections of $E$ is denoted by $\Gamma(E)$. Specifically, we denote $\Gamma(\wedge^{k}TM)$ as $\mathfrak{X}^{k}(M)$ and $\Gamma(\wedge^{k}T^{*}M)$ as $\Omega^{k}(M)$.
 The main references are \cite{MR2954043,MR2906391,MR4221224,MR0575344}.

\subsection{Geometric modules}
In algebraic geometry, K\"ahler differentials are used to develop differential theory over commutative algebras or schemes. For a commutative algebra $C$ over a base field $\mathbbm{k}$, recall that the module of \textit{K\"ahler differentials} of $C$ consists of a $C$-module $\Omega_{\mathbbm{k}}^{1}(C)$ and a $\mathbbm{k}$-derivation $d_{\text{alg}}:C\to \Omega_{\mathbbm{k}}^{1}(C)$ (called the \textit{universal derivation}) such that for any $C$-module $W$ and a $\mathbbm{k}$-derivation $D:C\to W$, there is a unique $C$-module homomorpism $\varphi: \Omega_{\mathbbm{k}}^{1}(C)\to W$ satisfying $\varphi d_{\text{alg}}=D$. For any $k\in \mathbb{N}$, we denote $\wedge_{C}^{k}\Omega_{\mathbbm{k}}^{1}(C)$ by $\Omega_{\mathbbm{k}}^{k}(C)$. The elements of $\Omega_{\mathbbm{k}}^{k}(C)$ are referred to as \textit{K\"ahler $k$-forms}. For any $C$-module $W$, a skew-symmetric $\mathbbm{k}$-multilinear map
$$D: \underbrace{C\times\cdots \times C}_{k \text{ copies}}\to W$$
is said to be a \textit{skew-symmetric $k$-derivation} with values in $W$, if $D$ is a $\mathbbm{k}$-derivation in each argument. The set of all skew-symmetric $k$-derivation with values in $W$ is denoted by $\mathfrak{X}^{k}(W)$. Clearly, $\mathfrak{X}^{0}(W)=W$ and $\mathfrak{X}^{1}(W)=\text{Der}_{\mathbbm{k}}(C,W)$. When $C=C^{\infty}(M)$ and $W=C$, $\mathfrak{X}^{k}(C^{\infty}(M))$ can be identified with $\mathfrak{X}^{k}(M)$ canonically, since we work in the context of real manifolds. An important property of higher K\"ahler differentials is the following. 
\begin{lemma}\label{highKahler-lem}
Let $C$ be a commutative $\mathbbm{k}$-algebra. For any $C$-module $W$ and $D\in \mathfrak{X}^{k}(W)$, there is a unique $C$-module homomorphism $\varphi$ such that the following diagram commutes:
$$\begin{tikzcd}
C^{k} \arrow[rr, "\wedge^{k}d_{\text{alg}}"] \arrow[rd, "D"'] &        & \Omega_{\mathbbm{k}}^{k}(C) \arrow[ld, "\varphi", dashed] \\
                                                                         & W&                                            
\end{tikzcd}$$
Here $\wedge^{k}d_{\text{alg}}:C^{k}\to \Omega_{\mathbbm{k}}^{k}(C),(a_{1},...,a_{k})\mapsto d_{\text{alg}}a_{1}\wedge\cdots\wedge d_{\text{alg}}a_{k}.$ 
\end{lemma}
However, when $C=C^{\infty}(M)$, $\Omega^{k}(M)$ does not possess the aforementioned universal property in the category of $C^{\infty}(M)$-modules. This difficulty can be overcome if one restricts oneself to the category of geometric $C^{\infty}(M)$-modules in the sense of Krasil'shchik-Lychagin-Vinogradov \cite{MR0861121}.

\begin{definition}
[See \cite{MR0861121,MR4221224}.]A $C^{\infty}(M)$-module\label{geo-module-def} $W$ is said to be \textit{geometric} if
$$\bigcap\limits_{p\in M}\mathfrak{m}_{p}W=0,$$
where $\mathfrak{m}_{p}$ is the maximal ideal of $C^{\infty}(M)$ consisting of the functions vanishing at $p\in M$.
\end{definition}
It is easy to see that any submodule of a geometric module is still geometric and for any $\mathbb{R}$-linear space $V$ and geometric $C^{\infty}(M)$-module $W$, $\text{Hom}_{\mathbb{R}}(V,W)$ is geometric. In particular, $\text{End}_{\mathbb{R}}W$ is geometric. An important example of geometric module is that for any vector bundle over $M$, say $E$, $\Gamma(E)$ is a geometric $C^{\infty}(M)$-module.
\begin{remark}\label{geomodule-rmk-nonprojfg}
However, this is not the only way in which geometric modules arise. For example, the infinite direct sum $\oplus_{i=1}^{\infty}C^{\infty}(M)$ is a geometric module, which is not finitely generated as a $C^{\infty}(M)$-module. In contrast, according to Serre-Swan theorem, for any vector bundle $E$ over $M$, $\Gamma(E)$ is a finitely generated projective $C^{\infty}(M)$-module. Clearly, any projective $C^{\infty}(M)$-module is geometric. In general, there are nonprojective geometric modules and nongeometric ones, see for example \cite[12.46]{MR4221224}.
\end{remark}
One can produce a geometric module from a given $C^{\infty}(M)$-module in a natural way.
\begin{example}\label{eg-geomfunctor}
[Geometrization, \cite{MR4221224}] Let $W$ be a $C^{\infty}(M)$-module. Then
$$\mathscr{G}(W):=W/\bigcap\limits_{p\in M}\mathfrak{m}_{p}W$$ 
gives a geometric $C^{\infty}(M)$-module, which is called the \textit{geometrization} of $W$. 

If $f:W\to W^{\prime}$ is a $C^{\infty}(M)$-module homomorphism, then $f$ induces a natural homomorphism $\mathscr{G}(f):\mathscr{G}(W)\to\mathscr{G}(W^{\prime})$. In this way we obtain a functor $\mathscr{G}$ from the category of $C^{\infty}(M)$-modules to the category of geometric $C^{\infty}(M)$-modules. This functor $\mathscr{G}$ is referred to as the \textit{geometrization functor}. Clearly, a $C^{\infty}(M)$-module $W$ is geometric if and only if $W\cong \mathscr{G}(W)$.
\end{example}
\begin{remark}\label{geo-Kahler-smooth1-rmk}
Unlike the case of multi-derivations, in general, the module of K\"ahler differentials $\Omega_{\mathbb{R}}^{1}(C^{\infty}(M))$ is not geometric. However, the module of smooth 1-forms $\Omega^{1}(M)$ is, in fact, the geometrization of $\Omega_{\mathbb{R}}^{1}(C^{\infty}(M))$, see \cite[14.11]{MR4221224}. Thus, it can be said that the module of “geometric differentials” is the geometrization of the module of “algebraic differentials”. This relationship also holds for higher differential forms, see Remark \ref{geome-alggeo-diffmo-rmk}.
\end{remark}
Later, in the construction of the complex (\ref{intro-poiss-complex}), we will apply the fact that the tensor product of any geometric module with the section module of a vector bundle is also geometric. This follows from the next more general result.
\begin{proposition}\label{tensor-geom-fgprojprop}
Let $W$ be a geometric $C^{\infty}(M)$-module and $P$ be a finitely generated projective $C^{\infty}(M)$-module. Then $W\otimes_{C^{\infty}(M)}P$ is also geometric. In particular, for any vector bundle $E$ over $M$, $W\otimes_{C^{\infty}(M)}\Gamma(E)$ is a geometric $C^{\infty}(M)$-module.
\end{proposition}
\begin{proof}
In general, for any ring $C$, any right $C$-module $X$ with a family of submodules $\{X_{i}\mid i\in\Lambda\}$, and any finitely generated projective left $C$-module $P$, the following holds:

$$\left(\bigcap\limits_{i\in \Lambda}X_{i}\right)\otimes_{C}P=\bigcap\limits_{i\in\Lambda}\left(X_{i}\otimes_{C}P\right)\subseteq X\otimes_{C}P.$$
Therefore, by setting $C=C^{\infty}(M), X=W,X_{i}=\mathfrak{m}_{p}W$, the conclusion follows.
\end{proof}


\par For any $k\in\mathbb{N}$, one has a canonical pairing 
\begin{align*}
\mathfrak{X}^{k}(M)\otimes_{C^{\infty}(M)}\Omega^{k}(M) &\to C^{\infty}(M)\\
X\otimes \omega &\mapsto  \hat{X}(\omega),
\end{align*}
where $\hat{X}:\Omega^{k}(M)\to C^{\infty}(M)$ is the unique $C^{\infty}(M)$-module homomorphism such that
$$\hat{X}(df_{1}\wedge\cdots\wedge df_{k})=X(df_{1},...,df_{k})$$
for all $f_{1},...,f_{k}\in C^{\infty}(M)$. The pairing induces a canonical isomorphism

\begin{equation}\label{duality-diff-cliso}
\mathfrak{X}^{k}(M)\to \text{Hom}_{C^{\infty}}(\Omega^{k}(M),C^{\infty}(M)),X\mapsto \hat{X}.
\end{equation}

When $k=1$, the duality (\ref{duality-diff-cliso}) is extended to the context of geometric modules, see \cite[14.7]{MR4221224}. This reflects the importance of geometric modules. Below is the corresponding result for higher differential forms, which can also be viewed as the geometric analogue of Lemma \ref{highKahler-lem}. Our proof is fundamentally an adaptation of the treatment for the case of 1-forms.
\begin{lemma}\label{unpro-highform-geolem}
For $k\in \mathbb{N}$, let $\wedge^{k}d$ be the natural skew-symmetric $k$-derivation induced by the exterior differential $d:C^{\infty}(M)\to\Omega^{1}(M)$, that is,$$\wedge^{k}d:C^{\infty}(M)\times\cdots \times C^{\infty}(M)\to \Omega^{k}(M),(f_{1},...,f_{k})\mapsto  df_{1}\wedge\cdots\wedge df_{k} .$$
Then for any geometric $C^{\infty}(M)$-module $W$ and $D\in \mathfrak{X}^{k}(W)$, there is a unique $C^{\infty}(M)$-module homomorphism  $\varphi:\Omega^{k}(M)\to W$ such that the following diagram commutes:
$$\begin{tikzcd}
C^{\infty}(M)^{k} \arrow[rr, "\wedge^{k}d"] \arrow[rd, "D"'] &        & \Omega^{k}(M) \arrow[ld, "\varphi", dashed] \\
                                                                         & W&                                            
\end{tikzcd}$$
In particular, one has a canonical isomorphism between $C^{\infty}(M)$-modules:
$$\Theta:\text{Hom}_{C^{\infty}(M)}(\Omega^{k}(M),W) \to\mathfrak{X}^{k}(W),\varphi \mapsto \varphi(\wedge^{k}d).$$
\end{lemma} 
\begin{proof}
Let $\Theta:\text{Hom}_{C^{\infty}(M)}(\Omega^{k}(M)),W)\to\mathfrak{X}^{k}(W)$ be the map defined by $\varphi\mapsto \varphi(\wedge^{k}d)$, which is clearly $C^{\infty}(M)$-linear. 
It suffices to construct the inverse map of $\Theta$. Since we  work in the context of real manifolds,  any $k$-form in $\Omega^{k}(M)$ can be expressed as $\sum_{i}g_{i0}dg_{i1}\wedge\cdots\wedge  dg_{ik} $ globally, see \cite[14.17]{MR4221224}. 
We claim that for fixed $p\in M$,
$$\Psi(D)(\sum_{i}g_{i0}dg_{i1}\wedge\cdots \wedge dg_{ik} )+\mathfrak{m}_{p}W=\sum\limits_{i}g_{i0}D(g_{i1},..., g_{ik})+\mathfrak{m}_{p}W $$
 does not depend on the choice of the representation $\sum_{i}g_{i0}dg_{i1}\wedge\cdots\wedge  dg_{ik}$. Once the assertion is proven, one can see that $\sum_{i}g_{i0}D(g_{i1},...,g_{ik})$ does not depend on the specific representation of $\sum_{i}g_{i0}dg_{i1}\wedge\cdots\wedge  dg_{ik}$, given that $W$ is geometric. Consequently, this allows one to define a well-defined map $\Psi:\mathfrak{X}^{k}(W)\to \text{Hom}_{C^{\infty}(M)}(\Omega^{k}(M),W)$ satisfying
 $$\Psi(D)(\sum_{i}g_{i0}dg_{i1}\wedge\cdots \wedge dg_{ik} )=\sum\limits_{i}g_{i0}D(g_{i1},..., g_{ik}).$$
Then it is straightforward to verify that $\Psi$ is indeed the inverse map of $\Theta$. Now we prove the aforementioned claim to complete the proof of this lemma.
 
 In fact, it suffices to show that when $\omega= \sum_{i} dh_{i1}\wedge\cdots\wedge dh_{ik} $ satisfies that $\omega_{p}=\sum_{i}(dh_{i1}\wedge\cdots\wedge dh_{ik})_{p}=0$, then $\sum_{i}D(h_{i1},..., h_{ik})\in\mathfrak{m}_{p}W$. This is because if one has
 
 $$\sum_{i}g_{i0}dg_{i1}\wedge\cdots\wedge dg_{ik}=\sum_{j}f_{j0}df_{j1}\wedge\cdots\wedge df_{jk}\in\Omega^{k}(M),$$
then $\sum_{i}g_{i0}(p)(dg_{i1})_{p}\wedge\cdots\wedge (dg_{ik})_{p}=\sum_{j}f_{j0}(p)(df_{j1})_{p}\wedge\cdots\wedge (df_{jk})_{p}\in T_{p}^{*}M$, which implies that 
 $$\sum_{i}(dg_{i0}(p)g_{i1})_{p}\wedge (dg_{i1})_{p}\wedge \cdots \wedge (dg_{ik})_{p}-\sum_{j}(df_{j0}(p)f_{j1})_{p}\wedge (df_{j2})_{p}\wedge\cdots\wedge (df_{jk})_{p}=0.$$
 
 Hence, the result follows once we show that $\omega_{p}=0$ implies that $\sum_{i}D(h_{i1},..., h_{ik})\in\mathfrak{m}_{p}W$. 
 
Assume that
 $\omega_{p}= \sum_{i}(dh_{i1})_{p}\wedge\cdots\wedge (dh_{ik})_{p}=0\in \wedge^{k}T_{p}^{*}M $. Then the action of $\omega_{p}$ on any $k$-vector in $\wedge^{k}T_{p}M$ is identically zero. Let $\{\gamma_{\alpha}|\alpha\in\Lambda_{1}\}$ be an $\mathbb{R}$-basis of $ \mathfrak{m}_{p}W$. Then we can extend this basis to an $\mathbb{R}$-basis of $W$, denoted by $\{\gamma_{\alpha}|\alpha\in\Lambda_{1}\}\cup  \{\gamma_{\beta}|\beta\in\Lambda_{2}\} $. For fixed $f_{1},...,f_{k}\in C^{\infty}(M)$, we write 
 $$D(f_{1},..., f_{k})=\sum_{\alpha\in \Lambda_{1}}X_{\alpha}(f_{1},...,f_{k})\gamma_{\alpha}+\sum_{\beta\in \Lambda_{2}}X_{\beta}(f_{1},...,f_{k})\gamma_{\beta},$$
 where $X_{\alpha}(f_{1},...,f_{k}),X_{\beta}(f_{1},...,f_{k})\in\mathbb{R}$ and all but finitely many coefficients of $\gamma_{\alpha},\gamma_{\beta}$ are zero. Clearly, $X_{\beta}$ is alternating $\mathbb{R}$-multilinear. For any $g_{1},g_{1}^{\prime},g_{2},...,g_{k}\in C^{\infty}(M)$ one has $X_{\beta}(g_{1}g_{1}^{\prime},g_{2},...,g_{k})=g_{1}(p)X_{\beta}(g_{1}^{\prime},g_{2},...,g_{k})+g_{1}^{\prime}(p)X_{\beta}(g_{1},g_{2},...,g_{k})$, since $g-g(p)\in\mathfrak{m}_{p}$ for any $g\in C^{\infty}(M)$.  Thus, $X_{\beta}:(C^{\infty}(M))^{k}\to C^{\infty}(M)$ defines a skew-symmetric $k$-derivation at $p$, for all $\beta\in\Lambda_{2}$. Hence, $X_{\beta}(\beta\in\Lambda_{2})$ can be viewed as an element in $\wedge^{k}T_{p}M$.

By assumption, the action of $\omega_{p}$ on any $k$-vector in $\wedge^{k}T_{p}M$ is zero. 
In particular, it follows that $\omega_{p}(X_{\beta})=0,\forall \beta\in\Lambda_{2}$. This can be rephrased in the following form:
 $$\sum\limits_{i}X_{\beta}(h_{i1},...,h_{ik})=0,\forall \beta\in\Lambda_{2}.$$
 Note that $\sum_{i}X_{\beta}(h_{i1},...,h_{ik})$ is the coordinate of $\sum_{i}D(h_{i1},...,h_{ik})+\mathfrak{m}_{p}W$ with respect to the given basis at $\gamma_{\beta}$. Therefore, it follows that $\sum_{i}D(h_{i1},..., h_{ik})\in\mathfrak{m}_{p}W$. 
\end{proof}
\begin{remark}\label{geome-alggeo-diffmo-rmk}
Consider the module of K\"ahler $k$-forms of $C^{\infty}(M)$, that is, $\Omega_{\mathbb{R}}^{k}(C^{\infty}(M))$. As pointed out in Remark \ref{geo-Kahler-smooth1-rmk}, in general, $\Omega^{k}(C^{\infty}(M))\neq \Omega^{k}(M)$. But from Lemma \ref{highKahler-lem} and Lemma \ref{unpro-highform-geolem}, we can immediately conclude that $\mathscr{G}(\Omega_{\mathbb{R}}^{k}(C^{\infty}(M)))\cong \Omega^{k}(M)$. Recall that $\mathscr{G}(\Omega_{\mathbb{R}}^{k}(C^{\infty}(M)))$ represents the geometrization of $\Omega_{\mathbb{R}}^{k}(C^{\infty}(M))$.
\end{remark}
Recall that for any vector bundle $E$ over $M$, $\Gamma(E\otimes \wedge^{k}TM)$, which is the space of $k$-vector fields on $M$ with values in $E$, can be identified with $\mathfrak{X}^{k}(\Gamma(E))$ canonically. This identification can be naturally extended to the category of geometric modules.
\begin{proposition}\label{tensorvec-mulder-prop}
Let $W$ be a geometric $C^{\infty}(M)$-module. Then the natural map
$$\Delta:W\otimes_{C^{\infty}(M)}\mathfrak{X}^{k}(M)\to \mathfrak{X}^{k}(W),w\otimes X\mapsto \Delta(w\otimes X), $$where $\Delta(w\otimes X)(f_{1},...,f_{k})=X(df_{1},...,df_{k})w$, is a $C^{\infty}(M)$-module isomorphism.
\end{proposition}
\begin{proof}
Since $\Omega^{k}(M)$ is finitely generated projective as a $C^{\infty}(M)$-module, there is a  canonical isomorphism
$$\Upsilon:  W\otimes_{C^{\infty}(M)}\text{Hom}_{C^{\infty}(M)}(\Omega^{k}(M),C^{\infty}(M)) \to \text{Hom}_{C^{\infty}(M)}(\Omega^{k}(M),W),$$
defined by $\Upsilon(w\otimes\varphi):\Omega^{k}(M)\to W,\eta\mapsto \varphi(\eta)w$ for all $\eta\in \Omega^{k}(M)$ and $C^{\infty}(M)$-linear map $\varphi:\Omega^{k}(M)\to C^{\infty}(M)$. Now Lemma \ref{unpro-highform-geolem} provides us with an isomorphism
$$\Theta:\text{Hom}_{C^{\infty}(M)}(\Omega^{k}(M),W)\to\mathfrak{X}^{k}(W),$$
which makes the following diagram commute
$$\begin{tikzcd}
{ W\otimes_{C^{\infty}(M)}\text{Hom}_{C^{\infty}(M)}(\Omega^{k}(M),C^{\infty}(M))} \arrow[rr, "\cong"] \arrow[d, "\Upsilon"'] &  & W\otimes_{C^{\infty}(M)}\mathfrak{X}^{k}(M) \arrow[d, "\Delta"] \\
{\text{Hom}_{C^{\infty}(M)}(\Omega^{k}(M),W)} \arrow[rr, "\Theta"]                                                            &  & \mathfrak{X}^{k}(W)                                            
\end{tikzcd}$$
It follows that $\Delta$ is a $C^{\infty}(M)$-module isomorphism.
\end{proof}

\subsection{Interior products}
Let $W$ be a geometric $C^{\infty}(M)$-module. By Proposition \ref{tensor-geom-fgprojprop}, one sees that $W\otimes_{C^{\infty}}\Omega^{k}(M)$ is also geometric. Consequently, Lemma \ref{unpro-highform-geolem} enables us to extend the construction of interior product induced by multivector fields.
\begin{definition}
[See also \cite{MR3395070}]For any $X\in\mathfrak{X}^{k}(W)$, the \textit{interior product} $\iota_{X}:\Omega^{q}(M)\to  W\otimes_{C^{\infty}(M)}\Omega^{q-k}(M)$ is a $C^{\infty}(M)$-module homomorphism of degree $-k$. It is defined as as follows: When $k\geq 1$, for $\omega=f_{0}df_{1}\wedge \cdots\wedge df_{q}\in\Omega^{q}(M)$,
\begin{equation}\label{def-interprodu}
\iota_{X}(\omega) = 
\begin{cases} 
0, & \text{if } q < k, \\
\sum_{\sigma \in S_{k,q-k}} \text{sgn}(\sigma) X(f_{\sigma(1)}, f_{\sigma(2)}, ..., f_{\sigma(k)}) f_0 \otimes df_{\sigma(k+1)} \wedge \cdots \wedge df_{\sigma(q)}, & \text{if } q \geq k,
\end{cases}
\end{equation}
where $S_{k,q-k}$ denotes the set of all $(k,q-k)$-shuffles, that is, $S_{k,q-k}=\{\sigma\in S_{q}\mid \sigma(1)<\cdots< \sigma(k),\sigma(k+1)<\cdots< \sigma(q)\}$. In the case when $k=0$, that is, $X$ is a smooth function on $M$, then $\iota_{X}:\Omega^{q}(M)\to\Omega^{q}(M)$ is defined by $\iota_{X}(\omega)=X\omega$.
\end{definition}
\begin{remark}
In some literatures, the interior product is also referred to as the \textit{contraction map}. If $W=C^{\infty}(M)$, one can identify $X\in\mathfrak{X}^{k}(C^{\infty}(M))$ with a smooth $k$-vector field on $M$, then the formula (\ref{def-interprodu}) can be rephrased as 
$$\iota_{X}(\omega)=\sum\limits_{\sigma\in S_{k,q-k}}\text{sgn}(\sigma)X(df_{\sigma(1)},df_{\sigma(2)},..., df_{\sigma(k)})f_{0} df_{\sigma(k+1)}\wedge\cdots\wedge df_{\sigma(q)}.$$
In this case, above definition boils down to the usual one. Recall also that the (generalized) \textit{Lie derivative} $\mathscr{L}_{X}$ of differential forms with respect to $X$ can be formulated by $\mathscr{L}_{X}=[\iota_{X},d]$, where $[-,-]$ denotes the graded commutator.
\end{remark}

We conclude this subsection by listing two key facts on the interior product that will be used later in the  proof of the main theorem.
\begin{proposition}
[\text{\cite[Proposition 3.4(2)]{MR2906391}}]For any $X\in\mathfrak{X}^{k}(W)$ and $Y\in\mathfrak{X}^{\ell}(M)$,\label{graded-commutat-iota-prop} $\iota_{X}\iota_{Y}=(-1)^{k\ell}(\text{id}_{W}\otimes \iota_{Y})\iota_{X}:\Omega^{q}(M)\to W\otimes_{C^{\infty}(M)}\Omega^{q-k-\ell}(M).$
\end{proposition}
The following can be viewed as the geometric analogue of \cite[Theorem 3.4]{MR4727073}.
\begin{proposition}
Suppose $M$ is a orientable smooth manifold of dimension $n$, and let $\mu$ be a volume form of $M$, then for any geometric module $W$ and integer $0\leq k\leq n$, the map $\star_{W}^{k}:\mathfrak{X}^{k}(W)\to W\otimes_{C^{\infty}(M)}\Omega^{n-k}(M),X\mapsto \iota_{X}\mu$ is a $C^{\infty}(M)$-module isomorphism\label{prop-starW-iso}.
\end{proposition}
\begin{proof}
In the case when $W=C^{\infty}(M)$, it is canonical that $\star^{k}: \mathfrak{X}^{k}(M)\to\Omega^{n-k}(M),X\mapsto \iota_{X}\mu$ is an isomorphism.  Consequently, the map
$$\text{id}_{W}\otimes \star^{k}:W\otimes_{C^{\infty}(M)}\mathfrak{X}^{k}(M)\to W\otimes_{C^{\infty}(M)}\Omega^{n-k}(M),w\otimes X\mapsto w\otimes\iota_{X}\mu$$
is also an isomorphism. It is straightforward to check that the following diagram is commutative:
$$\begin{tikzcd}
W\otimes_{C^{\infty}(M)}\mathfrak{X}^{k}(M) \arrow[rr, "\Delta"] \arrow[rd, "\text{id}_{W}\otimes \star^{k}"'] &                                         & \mathfrak{X}^{k}(W) \arrow[ld, "\star_{W}^{k}"] \\
 & W\otimes_{C^{\infty}(M)}\Omega^{n-k}(M) &                                                
\end{tikzcd}$$
where $\Delta$ is the isomorphism given in Proposition \ref{tensorvec-mulder-prop}. Thus, $\star_{W}^{k}$ is bijective.
\end{proof}

\section{Poisson geometric modules}\label{section 2}

In this section, we review some materials on Poisson Lie algebroids and Lie-Rinehart algebras that arise from Poisson algebras, as well as their representations. In the final part, we extend the correspondence, as established in \cite{MR1959580,MR2021638}, between the flat contravariant connections on a vector bundle over a Poisson manifold and the Poisson module structures on its section module to any geometric module.

\subsection{Poisson Lie algebroids and Lie-Rinehart algebras}

 Recall that a bivector field $\pi\in \mathfrak{X}^{2}(M)$ on a manifold $M$ is referred to as a \textit{Poisson structure}, if the bracket $\{-,-\}:C^{\infty}(M)\times C^{\infty}(M)\to C^{\infty}(M)$, defined by $\{f,g\}=\pi(df,dg)$ for $f,g\in C^{\infty}(M)$, is an $\mathbb{R}$-Lie bracket. When $\pi\in \mathfrak{X}^{2}(M)$ is a Poisson structure on $M$, the pair $(M,\pi)$ is called a \textit{Poisson manifold} and the bracket $\{-,-\}$ on $C^{\infty}(M)$ defined above is called the {\it Poisson bracket} on $M$. 

By axiomatizing the algebraic nature of the Poisson bracket on a Poisson manifold, one arrives at the concept of a Poisson algebra.

\begin{definition}
A commutative algebra $C$ over a field $\mathbbm{k}$, equipped with a bilinear map $\{-,-\}:C\times C\to C$, is called a \textit{Poisson algebra} if $(C,\{-,-\})$ is a $\mathbbm{k}$-Lie algebra and $\{-,-\}$ satisfies the Leibniz's rule, i.e., $\{ab,c\}=a\{b,c\}+\{a,c\}b$ for all $a,b,c\in C$.
\end{definition}

For a Poisson manifold $(M,\pi)$, there is a canonical bundle map $\pi^{\#}:T^{*}M\to TM$ defined by $(\pi^{\#})_{p}:T_{p}^{*}M\to T_{p}M, \omega_{p}\mapsto \iota_{\omega_{p}}(\pi_{p})$, where $\iota_{\omega_{p}}$ denotes the interior product by $\omega_{p}$. Then $\pi^{\#}$ induces a $C^{\infty}(M)$-linear map from $\Omega^{1}(M)$ to $\mathfrak{X}^{1}(M)$, which, by abuse of notation, is also denoted by  $\pi^{\#}$. For any $f\in C^{\infty}(M)$, the vector field $-\pi^{\#}(df)$ is referred to as the \textit{Hamiltonian vector field} of $f$, which is denoted by $X_{f}$. 

For a Poisson algebra $(C,\{-,-\})$, a $\mathbbm{k}$-derivation of the form $ \{-,a\}:C\to C$, defined by $b\mapsto \{b,a\}$, is referred to as the \textit{Hamiltonian derivation} of $a$. Suppose $f$ is a smooth function on a Poisson manifold $(M,\pi)$, then the derivation on the Poisson algebra $(C^{\infty}(M),\{-,-\})$ corresponding to the Hamiltonian vector field $X_{f}$ is precisely the Hamiltonian derivation $\{-,f\}$.

Poisson manifolds can be studied in the more general setting of Lie algebroids.

\begin{definition}
[\cite{MR216409}]A \textit{Lie algebroid} is a triple $(A,[-,-],\rho)$ consisting of a vector bundle $A$ on a manifold $M$, an $\mathbb{R}$-Lie algebra bracket $[-,-]$ on the space $\Gamma(A)$ and a bundle map $\rho:A\to TM$ (called the {\it anchor}), such that
\begin{itemize}
   \item[(1)] $\rho([X,Y])=[\rho(X),\rho(Y)]$;
   \item[(2)] $[X,fY]=f[X,Y]+(\rho(X)f)Y$
\end{itemize}
for all $X,Y\in\Gamma(A)$ and $f\in C^{\infty}(M)$.
\end{definition}
\begin{remark}
It is known that the condition (2) implies (1), see \cite{MR1077465}.
\end{remark}
For any smooth manifold $M$, the tangent bundle $TM$ naturally admits the structure of a Lie algebroid: One takes $[-,-]$ as the usual Lie bracket on vector fields and takes the anchor map as the identity map. The resulting Lie algebroid is called the \textit{tangent Lie algebroid} of $M$.
\begin{example}
[\cite{MR996653}]Let $(M,\pi)$ be a Poisson manifold \label{poiss-liealge-eg} and let $\pi^{\#}:T^{*}M \to TM$ be the bundle map induced by $\pi$. The Lie bracket $[-,-]$ on $\Omega^{1}(M)$ is defined by 
$$[\alpha,\beta]=\mathscr{L}_{\pi^{\#}(\alpha)}\beta-\mathscr{L}_{\pi^{\#}(\beta)}\alpha-d(\pi(\alpha,\beta)),$$
for all $\alpha,\beta\in \Omega^{1}(M)$, where $\mathscr{L}$ denotes the Lie derivative. Then $(T^{*}M,[-,-],\pi^{\#})$ is a Lie algebroid, called the \textit{Poisson Lie algebroid} or the \textit{contangent Lie algebroid} of $M$. 
\end{example}


One also has an algebraic analogue of the notion of Lie algebroid.
\begin{definition}
[\cite{MR0154906}]Let $C$ be a commutative $\mathbbm{k}$-algebra, and let $(L,[-,-])$ be a $\mathbbm{k}$-Lie algebra equipped with a left $C$-module structure. Let $\rho:L\to\text{Der}_{\mathbbm{k}}C$ be a Lie algebra homomorphism. The triple $(L,[-,-],\rho)$ is called a \textit{Lie-Rinehart algebra} over $\mathbbm{k}$ if
\begin{itemize}
    \item[(1)] $\rho(a\alpha)(b)=a\rho(\alpha)b$, for all $\alpha\in L,a,b\in C$, i.e., $\rho$ is $C$-linear;
    \item[(2)] $[\alpha,a\beta]=a[\alpha,\beta]+(\rho(\alpha)a)\beta$ for all $\alpha,\beta\in L,a\in C$.
\end{itemize}
For convenience, the pair $(C,L)$ is also referred to as a \textit{Lie-Rinehart algebra}.
\end{definition}

\begin{remark}
For any smooth manifold $M$, by indentifying $\mathfrak{X}(M)$ and $\text{Der}_{\mathbb{R}}C^{\infty}(M)$ canonically, one can easily see that any Lie algebroid $(A,[-,-],\rho)$ corresponds to a Lie-Rinehart algebra $(\Gamma(A),[-,-],\rho)$. For a Lie-Rinehart pair $(C,L)$ over $\mathbbm{k}$, $L$ is referred to as a \textit{$(\mathbbm{k},C)$-Lie algebra} by Rinehart \cite{MR0154906}. There are $(\mathbb{R},C^{\infty}(M))$-Lie algebras that do not arise from a Lie algebroid, since there are nongeometric $C^{\infty}(M)$-modules.
\end{remark}

Correspondingly, there is an algebraic version of Poisson Lie algebroid (Example \ref{poiss-liealge-eg}).
\begin{example}
[\text{\cite[Theorem 3.8]{MR1058984}}]Let $(C,\{-,-\})$ be a Poisson algebra over $\mathbbm{k}$.\label{eg-poissonalg-lirrine} Then one has a canonical $C$-module homomorphism $\rho:\Omega_{\mathbbm{k}}^{1}(C)\to \text{Der}_{\mathbbm{k}}C, d_{\text{alg}}a\mapsto \{a,-\}$. Additionally, there is a well-defined $\mathbbm{k}$-bilinear bracket $[-,-]:\Omega_{\mathbbm{k}}^{1}(C)\times \Omega_{\mathbbm{k}}^{1}(C)\to \Omega_{\mathbbm{k}}^{1}(C)$ such that
$$[ud_{\text{alg}}a,vd_{\text{alg}}b]=u\{a,v\}d_{\text{alg}}b+v\{u,b\}d_{\text{alg}}a+uvd_{\text{alg}}\{a,b\}$$
for all $u,v,a,b\in C$. In this case, $(C,\Omega_{\mathbbm{k}}^{1}(C))$ is a Lie-Rinehart algebra.
\end{example}
Given a Poisson manifold $(M,\pi)$, the Poisson Lie algebroid of $M$ produces a Lie-Rinehart algebra $(C^{\infty}(M),\Omega^{1}(M))$. Additionally, as a Poisson algebra, $(C^{\infty}(M),\{-,-\})$ defines another Lie-Rinehart algebra $(C^{\infty}(M),\Omega_{\mathbb{R}}^{1}(C^{\infty}(M)))$. In Remark \ref{geo-Kahler-smooth1-rmk}, we noted that $\Omega_{\mathbb{R}}^{1}(C^{\infty}(M))\neq \Omega^{1}(M) $ in general.  Consequently, $$(C^{\infty}(M),\Omega^{1}(M))\neq (C^{\infty}(M),\Omega_{\mathbb{R}}^{1}(C^{\infty}(M))).$$ We will see later  that for any geometric $C^{\infty}(M)$-module $W$, a flat $\Omega^{1}(M)$-connection on $W$ coincides with a flat $\Omega_{\mathbb{R}}^{1}(C^{\infty}(M))$-connection on $W$ (see Proposition \ref{geomodule-poisson-corresprop}).

Now we recall the concepts of connection and curvature of a Lie algebroid.
\begin{definition}
[\cite{MR896907,MR1675117}]Let $(A,[-,-],\rho)$ be a Lie algebroid over $M$ and let $E$ be a vector bundle on $M$. An $\mathbb{R}$-bilinear map $\nabla:\Gamma(A)\times \Gamma(E)\to\Gamma(E),(X,s)\mapsto \nabla_{X}s$ is called an \textit{$A$-connection} on $E$ if for any  $X\in\Gamma(A),s\in\Gamma(E)$, and $f\in C^{\infty}(M)$,
\begin{itemize}
    \item[(1)] $\nabla_{fX}(s)=f\nabla_{X}s$;
    \item[(2)] $\nabla_{X}(fs)=f\nabla_{X}(s)+(\rho(X)f)s$.
\end{itemize}
In this case, $ \Gamma(A)\times\Gamma(A)\to\text{End}_{\mathbb{R}}\Gamma(E),(X,Y)\mapsto \nabla_{X}\nabla_{Y}-\nabla_{Y}\nabla_{X}-\nabla_{[X,Y]}$ is called the \textit{curvature} of the $A$-connection $\nabla$. An $A$-connection is said to be \textit{flat} if its curvature is zero.
\end{definition}

\begin{remark}
When $\nabla$ is a flat $A$-connection on $E$, the pair $(E,\nabla)$ is also called a {\it representation} of the given Lie algebroid, see e.g., \cite{MR896907,MR1726784}.
\end{remark}

\begin{remark}
When $A$ is the tangent Lie algebroid of $M$, one recovers the classical notions of a (covariant) connection and curvature. For the Poisson Lie algebroid of a Poisson manifold $(M,\pi)$, a $T^{*}M$-connection on a vector bundle $E$ is known as the \textit{contravariant connection} on $E$ in the terminology of \cite{MR1818181} or the {\it contravariant derivative} on $E$ in \cite{MR1137387}.
\end{remark}
\begin{definition}
[\cite{MR1058984}]Let $(C,L )$ be a Lie-Rinehart algebra over $\mathbbm{k}$ and $W$ be a $C$-module. A $C$-linear map $\theta: L\to\text{End}_{\mathbbm{k}}W $ is said to be a {\it $L$-connection} on $W$ if
$$\theta(\alpha)(cw)=c\theta(\alpha)(w)+(\rho(\alpha)c)w$$
for all $c\in C, w\in W$, and $\alpha\in L$. As in the geometric case, the map $L\times L\to \text{End}_{\mathbbm{k}}W$, defined by $(\alpha,\beta)\mapsto \theta(\alpha)\theta(\beta)-\theta(\beta)\theta(\alpha)-\theta([\alpha,\beta])$, is called the \textit{curvature} of the $L$-connection $\theta$. And a $L$-connection $\theta$ is said to be \textit{flat}, if its curvature is zero.
\end{definition}
\begin{remark}When $W$ is equipped with a flat $L$-connection (in particular, $W$ is also a left $L$-module), $W$ is also referred to as a \textit{left $(C,L)$-module} by Huebschmann \cite{MR1058984}. Suppose $(C,L)$ is a Lie-Rinehart algebra, $W$ is a $C$-module equipped with a right $L$-module structure $[-,-]_{W}:W\times L\to L$ such that
$$
[cw,\alpha]_{W}=c[w,\alpha]_{W}-(\rho(\alpha)c)w=[w,c\alpha]_{W},\forall w\in W,c\in C,\alpha\in L.
$$
Then $W$ is referred to as a \textit{right $(C,L)$-module} \cite{MR1058984,MR1625610}. According to the definition, a $(C,L)$-module structure consists of a $C$-module structure and an $L$-module structure. However, when we consider a $(C,L)$-module structure on a $C$-module $W$, we always refer to the $L$-module structure that, together with the given $C$-module structure, makes $W$ a $(C,L)$-module. See also Remark \ref{lr-CLmodule-rmk}.
\end{remark}
\begin{remark}For any Lie algebroid $(A,[-,-],\rho)$ over $M$ and a vector bundle over $M$, consider the corresponding Lie-Rinehart algebra $(C^{\infty}(M),\Gamma(A))$. It is clear that a flat $A$-connection on $E$ coincides with a flat $\Gamma(A)$-connection on $\Gamma(E)$, or equivalently, a left $(C^{\infty}(M),\Gamma(A))$-module structure on $\Gamma(E)$.
\end{remark}

The representation theory and the (co)homology theory of Lie-Rinehart algebras can be studied via their universal enveloping algebras. We conclude this subsection by reviewing the definition of the universal enveloping algebra of a Lie-Rinehart algebra, which will be subsequently used to illustrate that the homology of the complex (\ref{intro-poiss-complex}) coincides with the geometric Poisson homology considered by Huebschmann.

\begin{definition}
[\cite{MR0154906,MR1058984}]Let $(C,L)$ be a Lie-Rinehart algebra over $\mathbbm{k}$. The \textit{universal enveloping algebra} $(\mathcal{U}(C,L),\iota_{C},\iota_{L})$ of $(C,L)$ is defined as a $\mathbbm{k}$-algebra $\mathcal{U}(C,L)$ together with a $\mathbbm{k}$-algebra homomorphism $\iota_{C}:C\to \mathcal{U}(C,L)$ and a $\mathbbm{k}$-Lie algebra homomorphism $\iota_{L}:L\to\mathcal{U}(C,L)$. The Lie bracket in $\mathcal{U}(C,L)$ is given by the commutator, and the following conditions are satisfied:
$$\iota_{C}(c)\iota_{L}(\alpha)=\iota_{L}(c\alpha),\  \iota_{L}(\alpha)\iota_{C}(c)-\iota_{C}(c)\iota_{L}(\alpha)=\iota_{C}(\rho(\alpha)c) ,\forall \alpha\in L,c\in C.$$
Moreover, for any $\mathbbm{k}$-algebra $\mathscr{B}$, which is also viewed as a $\mathbbm{k}$-Lie algebra, and for any pair of homomorphisms $\varepsilon_{C}:C\to \mathscr{B}$ (an algebra homomorphism) and $\varepsilon_{L}:L\to \mathscr{B}$ (a Lie algebra homomorphism) such that
$$\varepsilon_{C}(c)\varepsilon_{L}(\alpha)=\varepsilon_{L}(c\alpha),\  \varepsilon_{L}(\alpha)\varepsilon_{C}(c)-\varepsilon_{C}(c)\varepsilon_{L}(\alpha)=\varepsilon_{C}(\rho(\alpha)c) ,\forall \alpha\in L,c\in C,$$  
there exists a unique $\mathbbm{k}$-algebra homomorphism $\Phi:\mathcal{U}(C,L)\to \mathscr{B}$ such that $\Phi\iota_{C}=\varepsilon_{C}$ and $\Phi\iota_{L}=\varepsilon_{L}$, that is, the following diagram commutes:
$$\begin{tikzcd}
C \arrow[rr, "\iota_{C}"] \arrow[rrd, "\varepsilon_{C}"'] &  & {\mathcal{U}(C,L)} \arrow[d, "\Phi", dashed] &  & L \arrow[ll, "\iota_{L}"'] \arrow[lld, "\varepsilon_{L}"] \\
                                                   &  & \mathscr{B} &  &                                                   
\end{tikzcd}$$
\end{definition}
\begin{remark}\label{lr-CLmodule-rmk}
For the construction of the universal enveloping algebra of Lie-Rinehart algebras, see \cite{MR0154906,MR1058984}. Let $(C,L)$ be a Lie-Rinehart algebra over $\mathbbm{k}$ and let $W$ be a $C$-module. Then there is a natural one-to-one correspondence between left (right) $(C,L)$-module structures on $W$ and left (right) $\mathcal{U}(C,L)$-module structures compatible with the given $C$-module structure on $W$.  Here, the compatibility of the $\mathcal{U}(C,L)$-module structure with the given $C$-module structure on $W$ means that for any $c \in C$, the action of $\iota_{C}(c)$ on $W$ is consistent with the original action of $c$ on $W$. This is the reason for the distinction between left and right $(C,L)$-modules.
\end{remark}
\begin{remark}
When $L$ is a $\mathbbm{k}$-Lie algebra and $C=\mathbbm{k}$ with trivial $L$-action, then $\mathcal{U}(C,L)$ is the usual universal enveloping algebra of $L$. When $C=C^{\infty}(M)$ and $L=\mathfrak{X}(M)$, then $\mathcal{U}(C,L)$ is the \textit{algebra of (Grothendieck) differential operators} on $M$ \cite{MR1058984}. When $(C,\{-,-\})$ is a Poisson algebra and $L=\Omega_{\mathbbm{k}}^{1}(C)$, then the universal enveloping algebra of the Lie-Rinehart algebra in Example \ref{eg-poissonalg-lirrine} is the \textit{Poisson enveloping algebra} of $C$, introduced by Oh \cite{MR1683858}.
\end{remark}

\subsection{Flat contravariant connections and Poisson geometric modules}

Let $(M,\pi)$ be a Poisson manifold, fixed for the remainder of this paper. Then, we have the Poisson Lie algebroid $(T^{*}M,[-,-],\pi^{\#})$, the Poisson algebra $(C^{\infty}(M),\{-,-\})$ and two Lie-Rinehart algebras $(C^{\infty}(M),\Omega^{1}(M))$ and $(C^{\infty}(M),\Omega_{\mathbb{R}}^{1}(C^{\infty}(M)))$. In Remark \ref{lr-CLmodule-rmk}, we noted that for any $C^{\infty}(M)$-module $W$, there is a natural one-to-one correspondence between left (right) $\Omega^{1}(M)$-connections on $W$ and left (right) $\mathcal{U}(C^{\infty}(M),\Omega^{1}(M))$-module structures compatible with the given $C^{\infty}(M)$-module structure on $W$. Similarly, a left (right) $\Omega_{\mathbb{R}}^{1}(C^{\infty}(M))$-connection on $W$ corresponds to a left (right) $\mathcal{U}(C^{\infty}(M),\Omega_{\mathbb{R}}^{1}(C^{\infty}(M)))$-module structure compatible with compatible with the given $C^{\infty}(M)$-module structure on $W$, and vice versa. The aim of this subsection is to show that for any geometric $C^{\infty}(M)$-module $W$, there is a natural one-to-one correspondence between Poisson module structures on $W$, flat $\Omega^{1}(M)$-connections on $W$, and flat $\Omega_{\mathbb{R}}^{1}(C^{\infty}(M))$-connections on $W$ (Proposition \ref{geomodule-poisson-corresprop}).

In particular, when $W=\Gamma(E)$ is the module of global sections of a vector bundle $E$ over $M$, we see that a Poisson module structure on $\Gamma(E)$ can be identified with a flat $T^{*}M$-connection on $E$, which was proved in \cite {MR2021638,MR1959580}. Now, we begin with the definition of Poisson modules.
\begin{definition}
[\cite{MR1423640,MR1683858}]Let $(C,\{-,-\})$ be a Poisson algebra over $\mathbbm{k}$, and let $W$ be a $C$-module. Suppose $\{-,-\}_{W}:W\times C\to W$ is a bilinear map. The pair $(W,\{-,-\}_{W})$ is called a  right \textit{Poisson module} over $C$ provided:
\begin{itemize}
   \item[(1)] $(W,\{-,-\}_{W})$ is a right Lie module over the Lie algebra $(C,\{-,-\})$;
   \item[(2)] $\{wa,b\}_{W}=\{w,b\}_{W}a+w\{a,b\}$ for all $a,b\in C$ and $w\in W$;
   \item[(3)] $\{w,ab\}_{W}=\{w,a\}_{W}b+\{w,b\}_{W}a$ for all $a,b\in C$ and $w\in W$.
\end{itemize}
When $(W,\{-,-\}_{W})$ is a right Poisson module, the bilinear map $\{-,-\}_{W}$ is referred to as a right \textit{Poisson module structure} on $W$.
\end{definition}
\begin{remark}
The notion of left Poisson module is defined similarly. For any left Poisson module $(W,[-,-]_{W})$ over $C$,  $\{-,-\}_{W}:W\times C\to W,(w,c)\mapsto -[c,w]_{W}$ defines a right Poisson module structure on $W$. 
\end{remark}
\begin{remark}
In \cite{MR1423640}, a Poisson module over a Poisson algebra is referred to as a \textit{Poisson bimodule}. When $C=C^{\infty}(M)$ and $W=\Gamma(E)$ is the module of global sections of a vector bundle $E$ over $M$, then $E$, together with a Poisson module structure on $\Gamma(E)$, is also called a \textit{Poisson vector bundle} by Ginzburg \cite{MR1959580}.
\end{remark}
Given a Poisson algebra $(C,\{-,-\})$ and a $C$-module $W$, the definition of  K\"ahler differentials yields a one-to-one correspondence between right Poisson module structures on $W$ and flat $\Omega_{\mathbbm{k}}^{1}(C)$-connections on $W$. More precisely, for a right Poisson module structure on $W$, denoted by $\{-,-\}_{W}$, this structure defines a $\mathbbm{k}$-derivation $D:C\to\text{End}_{\mathbbm{k}}W$, given by $a\mapsto -\{-,a\}_{W}$. Let $d_{\text{alg}}:C\to\Omega^{1}(C)$ be the universal derivation, then there is a unique $C$-linear map $\theta:\Omega^{1}(C)\to\text{End}_{\mathbbm{k}}W$ such that the following diagram commutes.

$$\begin{tikzcd}
C \arrow[rd, "D"'] \arrow[rr, "d_{\text{alg}}"] &                           & \Omega_{\mathbbm{k}}^{1}(C) \arrow[ld, "\theta", dashed] \\
                                                & \text{End}_{\mathbbm{k}}W &                                                         
\end{tikzcd}$$

It is easy to check that $\theta$ is a flat $\Omega_{\mathbbm{k}}^{1}(C)$-connection on $W$. Conversely, any flat $\Omega_{\mathbbm{k}}^{1}(C)$-connection $\theta:\Omega_{\mathbbm{k}}^{1}(C)\to\text{End}_{\mathbbm{k}}W$ induces a right Poisson module structure on $W$ by setting 
$$\{w,a\}_{W}=-\theta(d_{\text{alg}}a)(w),\forall a\in C,w\in W.$$

\begin{proposition}\label{geomodule-poisson-corresprop}
Let $(M,\pi)$ be a Poisson manifold and $W$ be a geometric $C^{\infty}(M)$-module. Then there is a natural one-to-one correspondence between the following three sets:
\begin{itemize}
    \item[(1)] The set of all flat $\Omega^{1}(M)$-connections on $W$;
    \item[(2)] The set of all flat $\Omega_{\mathbb{R}}^{1}(C^{\infty}(M))$-connections on $W$;
    \item[(3)] The set of all right Poisson module structures on $W$.
\end{itemize}
We refer to a geometric $C^{\infty}(M)$-module equipped with a right Poisson module structure as a (right) \textit{Poisson geometric module}.
\end{proposition}
\begin{proof}
From the previous discussion, it suffices to prove that there is a one-to-one correspondence between flat $\Omega^{1}(M)$-connections on $W$ and flat $\Omega_{\mathbb{R}}^{1}(C^{\infty}(M))$-connections on $W$. Let $d_{\text{alg}}:C^{\infty}(M)\to\Omega_{\mathbb{R}}^{1}(C^{\infty}(M))$ be the universal derivation and let $d:C^{\infty}(M)\to \Omega^{1}(M)$ be the exterior differential. Then one has the canonical projection $\varkappa: \Omega_{\mathbb{R}}^{1}(C^{\infty}(M))\to \Omega^{1}(M)$, given by $g d_{\text{alg}}f\mapsto gdf $. Since $\Omega^{1}(M)$ is the geometrization of $\Omega_{\mathbb{R}}^{1}(C^{\infty}(M))$ (Remark \ref{geome-alggeo-diffmo-rmk}), we have $\text{Ker}\varkappa=\cap_{p\in M}\mathfrak{m}_{p}\Omega_{\mathbb{R}}^{1}(C^{\infty}(M))$.
\par Given a flat $\Omega^{1}(M)$-connection on $W$, say $\theta:\Omega^{1}(M)\to \text{End}_{\mathbb{R}}W$. By composing $\theta$ and $\varkappa$, we obtain a natural flat $\Omega_{\mathbb{R}}^{1}(C^{\infty}(M))$-connection on $W$. Conversely, for any flat $\Omega_{\mathbb{R}}^{1}(C^{\infty}(M))$-connection on $W$, say $\hat{\theta}:\Omega_{\mathbb{R}}^{1}(C^{\infty}(M))\to\text{End}_{\mathbb{R}}W$, one has $\cap_{p\in M}\mathfrak{m}_{p}\Omega_{\mathbb{R}}^{1}(C^{\infty}(M))\subseteq \text{Ker}\hat{\theta}$
because $\text{End}_{\mathbb{R}}W$ is a geometric $C^{\infty}(M)$-module. Thus $\hat{\theta}$ induces a unique flat $\Omega^{1}(M)$-connection $\theta:\Omega^{1}(M)\to \text{End}_{\mathbb{R}}W$ on $W$ such that $\theta\varkappa=\hat{\theta}$.
\end{proof}
\begin{remark}\label{rmk-poimodule-notE}
The class of Poisson geometric modules includes all smooth vector bundles with flat contravariant connections. However, not all Poisson geometric modules arise from vector bundles with flat contravariant connections, as there exist geometric modules that are not finitely projective. For example, consider a nonprojective geometric $C^{\infty}(M)$-module $W$. If we regard $M$ as a Poisson manifold with the trivial Poisson structure, then $W$ can be naturally viewed as a Poisson module with the trivial action.
\end{remark}
\begin{corollary}
[\cite{MR1959580,MR2021638}]Let $E$ be a vector bundle over $M$\label{cor-corres-module-flat}. There is a natural one-to-one correspondence between Poisson module structures on $\Gamma(E)$ over $C^{\infty}(M)$, left $(C^{\infty}(M),\Omega^{1}(M))$-module structures on $\Gamma(E)$, and flat $T^{*}M$-connections on $E$.
\end{corollary}
To relate the homology of the complex (\ref{intro-poiss-complex}) with Huebschmann's Poisson homology with coefficients later (see Remark \ref{rmk-coinc-Hueb-ho}), we need the following observation.
\begin{proposition}\label{prop-RLmodule-right}
Let $W$ be a geometric $C^{\infty}(M)$-module. Then there is a natural one-to-one correspondence between right Poisson module structures on $W$ and right $\mathcal{U}(C^{\infty}(M),\Omega^{1}(M))$-module structures compatible with the given $C^{\infty}(M)$-module structure on $W$.
\end{proposition}
\begin{proof}
For the sake of simplifying notation, denote $C^{\infty}(M)$ and $\Omega^{1}(M)$ as $C$ and $L$, respectively. From Lemma \ref{unpro-highform-geolem}, we see that any right Poisson module structure $ \{-,-\}_{W}:W\times C\to W$ corresponds to a $\mathbb{R}$-linear map $\zeta:L \to \text{End}_{\mathbb{R}}W$ such that $\zeta(df)=\{-,f\}_{W}$ for all $f\in C=C^{\infty}(M)$.  This induces a Lie algebra homomorphism $\psi_{L}=\zeta:L\to(\text{End}_{\mathbb{R}}W)^{op}$, where $(\text{End}_{\mathbb{R}}W)^{op}$ denotes the opposite algebra of $\text{End}_{\mathbb{R}}W$. One also has an $\mathbb{R}$-algebra homomorphism $\psi_{C}:C\to(\text{End}_{\mathbb{R}}W)^{op}$ defined by the left multiplication of elements in $C$ on  $W$. Applying the universal property of $(\mathcal{U}(C,L),\iota_{C},\iota_{L})$ to the triple $((\text{End}_{\mathbb{R}}W)^{op},\psi_{C},\psi_{L})$, one obtains the right $\mathcal{U}(C,L)$-module structure on $W$ corresponding to $\{-,-\}_{W}$.
\par Conversely, any right $\mathcal{U}(C,L)$-module structure on $W$ corresponds to an algebra homomorphism $\Psi:\mathcal{U}(C,L)\to (\text{End}_{\mathbb{R}}W)^{op}$. Suppose $\Psi$ is compatible with the scalar multiplication of $C$, that is, for any $f\in C$, $\Psi\iota_{C}(f)$ is the scalar multiplication of $f$ on $W$. Then $\Psi$ is a $C$-module homomorphism when $\mathcal{U}(C,L)$ is endowed with a left $C$-module structure via $\iota_{C}$. This leads to a $C$-module homomorphism $\Psi\iota_{L}:L\to (\text{End}_{\mathbb{R}}W)^{op}$. By setting $\{w,f\}_{W}=\Psi\iota_{L}(df)w$, for all $f\in C,w\in W$, it is straightforward to verify that $(W,\{-,-\}_{W})$ is a right Poisson module.
\end{proof}
\begin{remark}
For a general Lie-Rinehart algebra $(R,L)$, it is possible that there exist $R$-modules with a left $(R,L)$-module structure but do not admit a right $(R,L)$-module structure, see \cite[Theorem 1.1]{MR3366558} for example.
\end{remark}

\section{Poincar\'e duality between Poisson homologies and cohomologies}\label{section 3}
In this section, we first geometrize the constructions of the algebraic Poisson chain complexes with coefficients in Poisson modules to obtain the complex (\ref{intro-poiss-complex}), and subsequently, we geometrize the constructions of twisted Poisson geometric modules introduced in \cite{MR3395070}. Thanks to the constructions of twisted Poisson geometric modules, we prove that for any orientable Poisson manifold $M$, there is an explicit isomorphism between the Poisson cochain complex of $M$ with coefficients in any Poisson geometric module and the Poisson chain complex of $M$ with coefficients in the corresponding twisted Poisson geometric module. Consequently, a version of twisted Poincar\'e duality is established between the Poisson homologies and the Poisson cohomologies of $M$ with coefficients in any Poisson geometric module. In the case when the Poisson geometric module is induced by a vector bundle with a flat $T^{*}M$-connection, the aforementioned result generalizes \cite[Corollary 4.6]{MR1726784} and \cite[Theorem 4.8]{MR1675117}.  Throughtout, we fix a Poisson manifold $(M,\pi)$ and a (right) Poisson geometric module $(W,\{-,-\}_{W})$.
\subsection{Poisson (co)chain complexes with coefficients in Poisson geometric modules}
In this subsection, we study the Poisson chain complexes and cochain complexes with coefficients in Poisson geometric modules. 
\par In algebraic context, for any $\mathbbm{k}$-Poisson algebra $(C,\{-,-\})$ and a right Poisson module $(P,\{-,-\}_{P})$ over $C$, there is a canonical chain complex \cite{MR1647186,MR2275449,maszczyk2006}
\begin{equation}\label{alg-Poisson-complex}
\begin{tikzcd}[column sep=small]
\cdots \arrow[r] & P\otimes_{C}\Omega_{\mathbbm{k}}^{q}(C) \arrow[r, "\partial_{q}^{\text{alg}}"] & P\otimes_{C}\Omega_{\mathbbm{k}}^{q-1}(C) \arrow[r] & \cdots \arrow[r, "\partial_{2}^{\text{alg}}"] & P\otimes_{C}\Omega_{\mathbbm{k}}^{1}(C) \arrow[r, "\partial_{1}^{\text{alg}}"] & P \arrow[r] & 0,
\end{tikzcd}
\end{equation}
where $\partial_{r}^{\text{alg}}: P\otimes_{C}\Omega_{\mathbbm{k}}^{r}(C)\to P\otimes_{C}\Omega_{\mathbbm{k}}^{r-1}(C)$ is defined as:
\begin{align*}
\partial_{r}^{\text{alg}}(x\otimes& a_{0}d_{\text{alg}}a_{1}\wedge\cdots\wedge d_{\text{alg}}a_{r})=\sum\limits_{i=1}^{r}(-1)^{i-1}\{xa_{0},a_{i}\}_{P}\otimes d_{\text{alg}}a_{1}\wedge\cdots\widehat{d_{\text{alg}}a_{i}}\cdots\wedge d_{\text{alg}}a_{r}\\
+&\sum\limits_{1\leq i< j\leq r}(-1)^{i+j}x\otimes a_{0}d_{\text{alg}}\{a_{i},a_{j}\}\wedge d_{\text{alg}}a_{1}\wedge\cdots\widehat{d_{\text{alg}}a_{i}}\cdots\widehat{d_{\text{alg}}a_{j}}\cdots\wedge d_{\text{alg}}a_{r}.\end{align*}
The complex (\ref{alg-Poisson-complex}) is called the (algebraic) {\it Poisson chain complex} of $C$ with coefficients in $P$. Next we geometrize the complex (\ref{alg-Poisson-complex}) in the setting of Poisson manifolds. To simplify the notation, we will denote $C^{\infty}(M)$ as $C$ from now.
\begin{proposition} Let $W$ be a geometric Poisson module. 
There is a well-defined chain complex\label{prop-poisson-chain}
\begin{equation}\label{geo-poisson-complex}
\begin{tikzcd}[column sep=small]
\cdots \arrow[r] & W\otimes_{C}\Omega^{q}(M) \arrow[r, "\partial_{q}"] & W\otimes_{C}\Omega^{q-1}(M) \arrow[r] & \cdots \arrow[r, "\partial_{2}"] & W\otimes_{C}\Omega^{1}(M) \arrow[r, "\partial_{1}"] & W \arrow[r] & 0,
\end{tikzcd}
\end{equation}
where $\partial_{r}: W\otimes_{C}\Omega^{r}(M)\to W\otimes_{C}\Omega^{r-1}(M)$ is defined as:
\begin{align*}
\partial_{r}(w\otimes f_{0}df_{1}\wedge\cdots\wedge df_{r})=&\sum\limits_{i=1}^{r}(-1)^{i-1}\{wf_{0},f_{i}\}_{W}\otimes df_{1}\wedge\cdots\widehat{df_{i}}\cdots\wedge df_{r}\\
&+\sum\limits_{1\leq i< j\leq r}(-1)^{i+j}w\otimes f_{0}d\{f_{i},f_{j}\}\wedge df_{1}\wedge\cdots\widehat{df_{i}}\cdots\widehat{df_{j}}\cdots\wedge df_{r}.
\end{align*}
\end{proposition}
We refer to the complex (\ref{geo-poisson-complex}) as the (geometric) \textit{Poisson chain complex} of $M$ with coefficients in the geometric Poisson module $W$. The $r$-th homology of this complex is called the $r$-th \textit{Poisson homology} of $M$ with coefficients in $W$, denoted by $\text{HP}_{r}(M,W)$. In particular, when $W=C^{\infty}(M)$ (with the obvious Poisson module structure), this is precisely the canonical complex of $M$ introduced by Koszul \cite{MR0837203} and used in \cite{MR950556} to define the canonical homology of $M$.

\begin{proof}
It suffices to show that for each positive integer $r$, the map $\partial_{r}$ is well-defined. Once this is proven, consider the canonical projection $\varkappa_{r}:\Omega_{\mathbb{R}}^{r}(C)\to\Omega^{r}(M)$. Then it is clear that the following diagram commutes:
$$
\begin{tikzcd}[column sep=small]
\cdots \arrow[r] & W\otimes_{C}\Omega_{\mathbb{R}}^{q}(C) \arrow[r, "\partial_{q}^{\text{alg}}"] \arrow[d, "\text{id}_{W}\otimes \varkappa_{q}"] & W\otimes_{C}\Omega_{\mathbb{R}}^{q-1}(C) \arrow[r] \arrow[d, "\text{id}_{W}\otimes \varkappa_{q-1}"] & \cdots \arrow[r, "\partial_{2}^{\text{alg}}"] & W\otimes_{C}\Omega_{\mathbb{R}}^{1}(C) \arrow[r, "\partial_{1}^{\text{alg}}"] \arrow[d, "\text{id}_{W}\otimes \varkappa_{1}"] & W \arrow[r] \arrow[d, "\text{id}_{W}"] & 0 \\
\cdots \arrow[r] & W\otimes_{C}\Omega^{q}(M) \arrow[r, "\partial_{q}"]                                                                           & W\otimes_{C}\Omega^{q-1}(M) \arrow[r]                                                                & \cdots \arrow[r, "\partial_{2}"]              & W\otimes_{C}\Omega^{1}(M) \arrow[r, "\partial_{1}"]                                                                           & W \arrow[r]                            & 0
\end{tikzcd}
$$
This implies that $\partial_{r}\partial_{r+1}=0$ is an immediate consequence of $\partial_{r}^{\text{alg}}\partial_{r+1}^{\text{alg}}=0$.
\par In fact, for each $w\in W$, there is a well-defined $\mathbb{R}$-linear map $\theta_{w}:\Omega^{r}(M)\to W\otimes_{C}\Omega^{r-1}(M)$ given by:
\begin{align*}
\theta_{w}(f_{0}df_{1}\wedge\cdots\wedge df_{r})=&\sum\limits_{i=1}^{r}(-1)^{i-1}\{wf_{0},f_{i}\}_{W}\otimes df_{1}\wedge\cdots\widehat{df_{i}}\cdots\wedge df_{r}\\
&+\sum\limits_{1\leq i< j\leq r}(-1)^{i+j}w\otimes f_{0}d\{f_{i},f_{j}\}\wedge df_{1}\wedge\cdots\widehat{df_{i}}\cdots\widehat{df_{j}}\cdots\wedge df_{r}.
\end{align*}
To see this, let $\partial_{r}^{\text{can}}:\Omega^{r}(M)\to \Omega^{r-1}(M)$ be the differential of the canonical complex of $M$, that is, $\partial^{\text{can}}=[\iota_{\pi},d]$. Then one obtains an $\mathbb{R}$-linear map $w\otimes \partial_{r}^{\text{can}}:\Omega^{r}(M)\to W\otimes_{C}\Omega^{r-1}(M)$ such that
\begin{align*}
(w\otimes \partial_{r}^{\text{can}})(f_{0}df_{1}\wedge\cdots\wedge df_{r})=&\sum\limits_{i=1}^{q}(-1)^{i-1}w \{f_{0},f_{i}\}\otimes df_{1}\wedge\cdots \widehat{df_{i}}\cdots\wedge df_{r}\\
&+\sum\limits_{1\leq i< j\leq r}(-1)^{i+j }w \otimes f_{0}d\{f_{i},f_{j}\}\wedge df_{1}\wedge\cdots \widehat{df_{i}}\cdots \widehat{df_{j}}\cdots\wedge df_{r}.
\end{align*}
Note that $X=\{w,-\}_{W}\in\mathfrak{X}^{1}(W)$, one also has the interior product $\iota_{X}$ induced by $X$. Set $$\theta_{w}=w\otimes \partial_{r}^{\text{can}}+\iota_{X}.$$
Then
\begin{align*}
\theta_{w}(f_{0}df_{1}\wedge\cdots\wedge df_{r})=&\sum\limits_{i=1}^{q}(-1)^{i-1}w \{f_{0},f_{i}\}\otimes df_{1}\wedge\cdots \widehat{df_{i}}\cdots\wedge df_{r}\\
&+\sum\limits_{1\leq i< j\leq r}(-1)^{i+j }w \otimes f_{0}d\{f_{i},f_{j}\}\wedge df_{1}\wedge\cdots \widehat{df_{i}}\cdots \widehat{df_{j}}\cdots\wedge df_{r}\\
 &+\sum\limits_{i=1}^{r}(-1)^{i-1}\{w,f_{i}\}_{W}f_{0}\otimes df_{1}\wedge\cdots\widehat{df_{i}}\cdots\wedge df_{r}\\
 =&\sum\limits_{i=1}^{r}(-1)^{i-1}(w\{f_{0},f_{i}\}+\{w,f_{i}\}_{W}f_{0})\otimes df_{1}\wedge\cdots\widehat{df_{i}}\cdots\wedge df_{r}\\
&+\sum\limits_{1\leq i< j\leq r}(-1)^{i+j}w\otimes f_{0}d\{f_{i},f_{j}\}\wedge df_{1}\wedge\cdots\widehat{df_{i}}\cdots\widehat{df_{j}}\cdots\wedge df_{r}\\
=&\sum\limits_{i=1}^{r}(-1)^{i-1}\{wf_{0},f_{i}\}_{W}\otimes df_{1}\wedge\cdots\widehat{df_{i}}\cdots\wedge df_{r}\\
&+\sum\limits_{1\leq i< j\leq r}(-1)^{i+j}w\otimes f_{0}d\{f_{i},f_{j}\}\wedge df_{1}\wedge\cdots\widehat{df_{i}}\cdots\widehat{df_{j}}\cdots\wedge df_{r}.
\end{align*}
Now one has an $\mathbb{R}$-bilinear map $W\times \Omega^{r}(M)\to \Omega^{r-1}(M),(w,\eta)\mapsto \theta_{w}(\eta)$, which induces the desired $\mathbb{R}$-linear map $\partial_{r}: W\otimes_{C}\Omega^{r}(M)\to W\otimes_{C}\Omega^{r-1}(M)$.
\end{proof}
\begin{remark}
The complex (\ref{geo-poisson-complex}) cannot be derived  by directly applying the tensor functor $W\otimes_{C}-$ to the canonical complex of $M$, because the differential of the canonical complex is not $C$-linear.
\end{remark}
\begin{remark}\label{rmk-coinc-Hueb-ho}
It is possible to apply Proposition \ref{prop-RLmodule-right} to construct the complex (\ref{geo-poisson-complex}) from a complex constructed by Rinehart in the setting of Lie-Rinehart algebras. For any Lie-Rinehart algebra $(R,L)$, Rinehart \cite[Eq. (4.1)]{MR0154906} constructed the following complex:
\begin{equation}\label{Rinehart-complex}
\begin{aligned}
\cdots \longrightarrow \, {\mathcal{U}(R,L)\otimes_{R}\wedge_{R}^{q}L} \,\stackrel{\mathfrak{d}_{q}}{\longrightarrow}\, & {\mathcal{U}(R,L)\otimes_{R}\wedge_{R}^{q-1}L} \stackrel{\mathfrak{d}_{q-1}} {\longrightarrow} \cdots\\
\longrightarrow \, & {\mathcal{U}(R,L)\otimes_{R}L} \stackrel{\mathfrak{d}_{1}}{\longrightarrow}  {\mathcal{U}(R,L)\otimes_{R}R} \longrightarrow 0,
\end{aligned}
\end{equation}
where $\mathfrak{d}_{r}:\mathcal{U}(R,L)\otimes_{R}\wedge_{R}^{r}L\to \mathcal{U}(R,L)\otimes_{R}\wedge_{R}^{r-1}L$ is given by:
\begin{align*}
\mathfrak{d}_{r}(u\otimes \alpha_{1}\wedge\cdots\wedge \alpha_{r})=&\sum\limits_{i=1}^{r}(-1)^{i-1}u\iota_{L}(\alpha_{i})\otimes \alpha_{1}\wedge\cdots\widehat{\alpha_{1}}\cdots\wedge \alpha_{r}\\
&+\sum\limits_{1\leq i< j\leq r}(-1)^{i+j}u\otimes\alpha_{1}\wedge\cdots \widehat{\alpha_{i}}\cdots\widehat{\alpha_{j}}\cdots\wedge \alpha_{r},
\end{align*}
and $(\mathcal{U}(R,L),\iota_{R},\iota_{L}) $ is the universal enveloping algebra of $(R,L)$. When $R=C^{\infty}(M)$ and $L=\Omega^{1}(M)$, for any Poisson geometric module $(W,\{-,-\}_{W})$, by  Proposition \ref{prop-RLmodule-right}, $W$ has a unique right $\mathcal{U}(R,L)$-module structure such that $u\cdot \iota_{L}(df)=\{u,f\}_{W}$ and $u\cdot \iota_{C}(g)=ug$ for all $u\in \mathcal{U}(R,L)$ and $f,g\in R$. Applying the tensor product functor $W\otimes_{\mathcal{U}(R,L)}-$ to the complex (\ref{Rinehart-complex}), then we obtain the complex (\ref{geo-poisson-complex}). Rinehart also showed that when $L$ is a projective $R$-module, then the complex (\ref{Rinehart-complex}) defines a projective resolution of $R$ (as a left $\mathcal{U}(R,L)$-module), see \cite[Lemma 4.1]{MR0154906}. In particular, our notion of Poisson homology with coefficients in $(W,\{-,-\}_{W})$ coincides with the Huebschmann's Poisson homology of $M$ with coefficients in $W$ \cite{MR1058984} with the aforementioned right $\mathcal{U}(R,L)$-module structure.
\end{remark}
From Proposition \ref{geomodule-poisson-corresprop}, we have seen that for any vector bundle $E$ over $M$, flat $T^{*}M$-connections on $E$ are precisely the Poisson module structures on the geometric module $\Gamma(E)$. Therefore, Proposition \ref{prop-poisson-chain} enbales us to consider the Poisson chain complex with coefficients in any vector bundle with a flat $T^{*}M$-connection.
\begin{corollary}\label{cor-poissoncomplex-flat}
For any vector bundle $E$ over $M$ with a flat $T^{*}M$-connection $\nabla$, there is a well-defined chain complex
\begin{equation}\label{poisson-chain-flat}
\begin{tikzcd}
\cdots \arrow[r] & \Gamma(E\otimes \wedge^{q}T^{*}M)   \arrow[r] & \cdots \arrow[r, "\partial_{2}"] & \Gamma(E\otimes T^{*}M) \arrow[r, "\partial_{1}"] & \Gamma(E) \arrow[r] & 0
\end{tikzcd}
\end{equation}
where $\partial_{r}: \Gamma(E\otimes \wedge^{r}T^{*}M) \to \Gamma(E\otimes \wedge^{r-1}T^{*}M)$ is given by the formula
\begin{align*}
\partial_{r}(w\otimes f_{0}df_{1}\wedge\cdots\wedge df_{r})=&\sum\limits_{i=1}^{r}(-1)^{i-1}\nabla_{-df_{i}}(f_{0}w) \otimes df_{1}\wedge\cdots\widehat{df_{i}}\cdots\wedge df_{r}\\
&+\sum\limits_{1\leq i< j\leq r}(-1)^{i+j}w\otimes f_{0}d\{f_{i},f_{j}\}\wedge df_{1}\wedge\cdots\widehat{df_{i}}\cdots\widehat{df_{j}}\cdots\wedge df_{r}.\end{align*}
The complex (\ref{poisson-chain-flat}) is referred to as the \textit{Poisson chain complex} of $M$ with coefficients in $(E,\nabla)$. The $r$-th homology of this complex is called the $r$-th \textit{Poisson homology} of $M$ with coefficients in $(E,\nabla)$, denoted by $\text{HP}_{r}(M,\nabla)$.
\end{corollary}
The notions of Poisson cohomology and Poisson cochain complex can be defined algebraically, since it had been noted by Huebschmann that the algebraic Poisson cochain complex and the geometric Poisson cochain complex are essentially the same, see \cite[Theorem 3.12.13]{MR1058984}.
\begin{definition}
[\cite{MR0501133,MR1058984}]For $C=C^{\infty}(M)$, there is a cochain complex
\begin{equation}\label{poisson-cochain-poissonmodule}
\begin{tikzcd}
0 \arrow[r] & W \arrow[r, "\delta^{0}"] & \mathfrak{X}^{1}(W) \arrow[r, "\delta^{1}"] & \mathfrak{X}^{2}(W) \arrow[r, "\delta^{2}"] & \cdots,
\end{tikzcd}
\end{equation}
where $\delta^{r}:\mathfrak{X}^{r}(W)\to\mathfrak{X}^{r+1}(W)$ is defined as
\begin{align*}
\delta^{r}(X)(f_{1},f_{2},...,f_{r+1})=&\sum\limits_{i=1}^{r+1}(-1)^{i}\{X(f_{1},...,\hat{f_{i}},...,f_{r+1}),f_{i}\}_{W}\\
+&\sum\limits_{1\leq i< j\leq r+1}(-1)^{i+j}X(\{f_{i},f_{j}\},f_{1},...,\hat{f_{i}},...,\hat{f_{j}},...,f_{r+1}).
\end{align*}
The complex (\ref{poisson-cochain-poissonmodule}) is called the \textit{Poisson cochain complex} of $M$ with coefficients in $W$. The $r$-th cohomology of this complex is called the $r$-th \textit{ Poisson cohomology} of $M$ with coefficients in $W$, denoted by $\text{HP}^{r}(M,W)$. 
\end{definition}

When $E$ is a vector bundle over $M$ with a flat $T^{*}M$-connection $\nabla$, then it is clear that the  Poisson cochain complex with coefficients in the corresponding geometric Poisson module $\Gamma(E)$  is isomorphic to the cochain complex of the Poisson Lie algebroid $T^{*}M$ with coefficients in the representation $(E,\nabla)$, see for example \cite{MR1726784}. Hence $\text{HP}^{r}(M,\Gamma(E))$ is precisely the cohomology of $T^{*}M$ with coefficients in $(E,\nabla)$. In this case, $\text{HP}^{r}(M,\Gamma(E))$ is denoted by $\text{HP}^{r}(M,\nabla)$.

The elements in $\text{Ker}\delta^{1}$ are referred to as \textit{Poisson derivations}, while the elements in $\text{Im}\delta^{0}$ are called \textit{Hamiltonian derivations}, which take the form $\{w,-\}_{W}$ for $w\in W$. When $W=C^{\infty}(M)$, these two concepts correspond to \textit{Poisson vector field} and \textit{Hamiltonian vector field}, respectively.

\subsection{Modular vector fields and twisted Poisson modules}
In this subsection, we first recall the concepts of modular vector field and modular class for Poisson manifolds.  Then, for any vector bundle $E$, a flat $T^{*}M$-connection $\nabla$ on $E$, and an arbitrary Poisson vector field, we introduce the construction of \textit{twisted flat $T^{*}M$-connection} of $\nabla$ by the given Poisson vector field. This construction can be viewed as a geometric analogue of the \textit{twisted Poisson module} by a Poisson derivation, as introduced in \cite{MR3395070}.
\begin{definition}
[\cite{MR0837203,MR1484598}]Let $(M,\pi)$ be an orientable Poisson manifold of dimension $n$ with a volume form $\mu$. The vector field $\phi_{\mu}$, defined by
$$\mathscr{L}_{X_{f}}\mu=\phi_{\mu}(f)\mu$$
for all $f\in C^{\infty}(M)$, is called the \textit{modular vector field} of $M$ with respect to $\mu$, where $\mathscr{L}_{X_{f}}$ is the Lie derivative of the Hamiltonian vector field $X_{f}$.
\end{definition}
It is well known that one can also defined the modular vector field with respect to $\mu$ as the unique vector field $\phi_{\mu}$ such that $\mathscr{L}_{\pi}\mu =\iota_{\phi_{\mu}}\mu$. Note that $\mathscr{L}_{\pi}=\iota_{\pi}d-d\iota_{\pi}$ is precisely the differential $\partial$ of the canonical complex of $M$, hence one has
\begin{equation}\label{modular-diff-voleq}
\partial\mu=\iota_{\phi_{\mu}}\mu.
\end{equation}
For affine smooth Poisson algebra with trivial canonical bundle, the equation (\ref{modular-diff-voleq}) can be verified by using the dual basis of the module of K\"ahler differentials of the given Poisson algebra, see \cite[Proposition 3.11]{MR4727073} (the reason why the conclusion there differs by a negative sign from (\ref{modular-diff-voleq}) is that the conventions of Hamiltonian derivation of a given function/element are different).

\par Recall also that  for an orientable Poisson manifold, the Poisson cohomology class of modular vector field is independent of the choice volume form, and this cohomology class is called the \textit{modular class} of the given Poisson manifold. This characteristic class is the obstruction to the existence of a volume form that is preserved by all Hamiltonian flows. Namely, the given Poisson manifold admits a volume form preserved by all Hamiltonian flows if and only if its modular class is zero. An orientable Poisson manifold is said to be {\it unimodular}, if its modular class is zero. For example, any simplectic manifold is unimodular since the Liouville form on it is preserved by all Hamiltonian flows.
\par In the following, we fix a right Poisson geometric module $(W,\{-,-\}_{W})$ over $C^{\infty}(M)$.
\begin{definition}
[\text{\cite[Proposition 2.7]{MR3395070}}]Let $\phi\in \mathfrak{X}^{1}(M)$ be a Poisson vector field. Define an $\mathbb{R}$-bilinear map $\{-,-\}_{W_{\phi}}:W\times C^{\infty}(M)\to W$ as
$$\{w,f\}_{W_{\phi}}=\{w,f\}_{W}+w\phi(f),\forall w\in W,f\in C^{\infty}(M),$$
then $(W,\{-,-\}_{W_{\phi}})$ is a Poisson geometric module, which is called the \textit{twisted Poisson geometric module} of $W$ by the Poisson vector field $\phi$, denoted by $W_{\phi}$.
\end{definition}

It is worth to pointing out that the operation of “twisted by a Poisson vector field” is invertible. More precisely, if we twist a Poisson geometric module $W$ by a Poisson vector field by $\phi$, then the “double twisted” Poisson module $(W_{\phi})_{-\phi}$ is indeed $W$ itself.

\begin{proposition}\label{prop-twist-flatc}
Let $E$ be a vector bundle over $M$ with a flat $T^{*}M$-connection $\nabla$. Then for any Poisson vector field $\phi$ over $M$, there is a flat $T^{*}M$-connection $\nabla^{\phi}$ defined by:
\begin{equation}\label{twisted-flat-connec}
\nabla_{df}^{\phi}(s)=\nabla_{df}(s)-\phi(f)s,\forall f\in C^{\infty}(M),s\in\Gamma(E).
\end{equation}
The flat $T^{*}M$-connection $\nabla^{\phi}$ given by (\ref{twisted-flat-connec}) is referred to as the {\it twisted flat $T^{*}M$-connection} of $\nabla$ by $\phi$.
\end{proposition}
\begin{proof}
By Corollary \ref{cor-corres-module-flat}, $\nabla$ induces a right Poisson module structure on $\Gamma(E)$ such that
$$\{s,f\}_{\Gamma(E)}=-\nabla_{df}(s),\forall s\in \Gamma(E), f\in C^{\infty}(M).$$
Twisting the above Poisson structure by the Poisson vector field $\phi$, then $\nabla_{df}^{\phi}$ is precisely the flat $T^{*}M$-connection corresponding to the new right Poisson module structure $\{-,-\}_{\Gamma(E)_{\phi}}$.
\end{proof}
Now we relate the construction above to the canonical flat $T^{*}M$-connection on the canonical line bundle  discovered in \cite[Eq. (20)]{MR1726784} and \cite[Eq. (18)]{MR1675117} respectively.

\begin{proposition}
Assume that $M$ is an orientable Poisson manifold of dimension $n$\label{ELW-twistedill-prop} with a  volume form $\mu$ and $\phi_{\mu}$ is the modular vector field of $M$ with respect to $\mu$. Consider the obvious flat $T^{*}M$-connection $\nabla_{0}$ on the canonical line bundle $\wedge^{n}T^{*}M$ defined by
$$(\nabla_{0})_{df}(\mu)=0,\forall f\in C^{\infty}(M).$$
Then the twisted flat $T^{*}M$-connection $\nabla_{0}^{\phi_{\mu}}$ of $\nabla$ by $\phi_{\mu}$ satisfies
$$(\nabla_{0}^{\phi_{\mu}})_{\omega}s=\omega\wedge d\iota_{\pi}(s),$$
for all $\omega\in\Omega^{1}(M),s\in \Omega^{n}(M)$.
\end{proposition}
\begin{proof}
Straightforward.
\end{proof}
\begin{remark}\label{rmk-flat-trilinebundle}
In fact, under the assumption of the proposition, one can easily check that any flat $T^{*}M$-connection on the trivial line bundle $\wedge^{n}T^{*}M$ arises from a twisted flat $T^{*}M$-connection of the obvious flat $T^{*}M$-connection (defined by a given volume form) with respect to some Poisson vector field, see \cite[Lemma 2.3]{MR3687261}. See also \cite{MR1625610,MR1675117} for a discussion of all flat connections on the canonical line bundle of a Lie-Rinehart algebra/Lie algebroid.
\end{remark}

\subsection{Twisted Poincar\'e duality between Poisson homologies and cohomologies}
The main result in this subsection is to establish an explicit isomorphism between the Poisson chain complex and the Poisson cochain complex with coefficients in any Poisson geometric module. This result can be seen as the geometric version of \cite[Theorem 3.12]{MR4727073}, which generalizes the duality theorems in \cite{MR1726784} and \cite{MR1675117} in the setting of orientable Poisson manifolds.


 To avoid  notational confusion, denote the differentials of the Poisson chain complex and Poisson cochain complex of $M$ with coefficients in $W$ by $\partial^{W}$ and $\delta_{W}$ respectively. To simplify the notation, let $C=C^{\infty}(M)$, $\partial=\partial^{C}$ and $\delta=\delta_{C}$.
\par Prior to presenting the statement and proof of the main theorem, we require two fundamental lemmas concerning the differentials of the Poisson chain complex and Poisson cochain complex with coefficients in any Poisson geometric module. 
\begin{lemma}
[\text{\cite[Lemma 3.6]{MR4727073}}]For any $X\in\mathfrak{X}^{k}(W)$\label{lem-inter-chaindiff}, one has
$$\iota_{X}\partial-(-1)^{k}\partial^{W}\iota_{X}=\iota_{(\delta_{W}X)}:\Omega^{\bullet}(M)\to W\otimes_{C}\Omega^{\bullet-k-1}(M).$$
\end{lemma}
\begin{proof}
This formula was obtained in the context of Poisson algebra. For the geometric case, the formula can be obtained using the following two methods.
\begin{itemize}
  \item[(1)]By \cite[14.17]{MR4221224}, it suffices to show that for any $\omega=f_{0}df_{1}\wedge\cdots\wedge df_{r}\in\Omega^{r}(M)$, one has $\iota_{X}\partial(\omega)-(-1)^{k}\partial^{W}\iota_{X}\omega=\iota_{(\delta_{W}X)}\omega$. This can be achieved by adapting the calculations for K\"ahler differentials in the Poisson algebra setting from \cite[Lemma 3.5]{MR4727073} step by step to the calculations for smooth differentials in the Poisson manifold setting.
  \item[(2)]Knowing that $\Omega^{k}(M)$ is the geometrization of $\Omega_{\mathbb{R}}^{k}(C)$, for all $k\in\mathbb{N}$, we can apply the canonical projection $\varkappa_{k}:\Omega_{\mathbb{R}}^{k}(C)\to\Omega^{k}(M)$ to transform the corresponding formula given in \cite[Lemma 3.6]{MR4727073} into the setting of Poisson manifolds. 
  \par To avoid confusion in notation, we denote the differential of the algebraic Poisson chain complex with coefficients in the Poisson module $W$ by $\partial_{\text{alg}}^{W}$, and the interior product induced by $X$ in the setting of commutative algebras (referred to as the contraction map in \cite[Definition 1.9]{MR4727073}) as $\iota_{X}^{\text{alg}}$. Then $\iota_{X}^{\text{alg}}\partial_{\text{alg}}-(-1)^{k}\partial_{\text{alg}}^{W}\iota_{X}^{\text{alg}}=\iota_{(\delta_{W}X)}^{\text{alg}}$, by \cite[Lemma 3.6]{MR4727073}. Clearly, for any $q\in \mathbb{N}$, one has the following commutative diagram:
  $$\begin{tikzcd}
\Omega_{\mathbb{R}}^{q}(C) \arrow[d, "\varkappa_{q}"'] \arrow[rr, "\iota_{X}^{\text{alg}}"] &  &  W\otimes_{C}\Omega_{\mathbb{R}}^{q-k}(C) \arrow[d, "\text{id}_{W}\otimes\varkappa_{q-k}"] \\
\Omega^{q}(M) \arrow[rr, "\iota_{X}"]                                          &  &  W\otimes_{C}\Omega^{q-k}(M)                                                 
\end{tikzcd}$$
Thus, we have
\begin{align*}
(\iota_{X}\partial-(-1)^{k}\partial^{W}\iota_{X})\varkappa_{q}=& \iota_{X}\varkappa_{q}\partial_{\text{alg}}-(-1)^{k}\partial^{W}(\text{id}_{W}\otimes \varkappa_{q-k})\iota_{X}^{\text{alg}}\\
=&(\text{id}_{W}\otimes \varkappa_{q-k})(\iota_{X}^{\text{alg}}\partial_{\text{alg}}-(-1)^{k}\partial_{\text{alg}}^{W}\iota_{X}^{\text{alg}})\\
=&(\text{id}_{W}\otimes \varkappa_{q-k})\iota_{(\delta_{W}X)}^{\text{alg}}\\
=&\iota_{(\delta_{W}X)}\varkappa_{q}.
\end{align*}
Therefore, since $\varkappa_{q}$ is surjective, one obtains that $\iota_{X}\partial-(-1)^{k}\partial^{W}\iota_{X}=\iota_{(\delta_{W}X)}$.
\end{itemize}
\end{proof}
\begin{lemma}
[\text{\cite[Lemma 3]{MR3969113}}]Let $\phi$ be a Poisson vector field on $M$\label{lem-twisted-diff}. Then
$$\delta_{W_{\phi}}=\delta_{W}-(\phi\wedge -) \text{ and \ }\partial^{W_{\phi}}=\partial^{W}+\text{id}_{W}\otimes \iota_{\phi}.$$
\end{lemma}
Suppose $M$ is an orientable Poisson manifold of dimension $n$ with a volume form $\mu$, $\nabla_{0}$ is the obvious flat $T^{*}M$-connection on the canonical line bundle $\wedge^{n}T^{*}M$ induced by $\mu$ in the sense of Proposition \ref{ELW-twistedill-prop}, and $\nabla_{0}^{-\phi_{\mu}}$ is the twisted flat $T^{*}M$-connection by the opposite modular vector field $-\phi_{\mu}$, then the differential of Poisson chain complex with coefficients in $(\wedge^{n}T^{*}M,\nabla^{-\phi_{\mu}})$ corresponds to the \textit{generating operator} of $\nabla_{0}$ considered by Xu, see \cite[Proposition 4.7]{MR1675117}. Thus   our notion of Poisson homology with coefficients in the twisted flat $T^{*}M$-connection $\nabla_{0}^{-\phi_{\mu}}$ coincides with the Lie algebroid homology with respect to the obvious flat $T^{*}M$-connection $\nabla_{0}$ on $\wedge^{n}T^{*}M$ in the Poisson case \cite{MR1675117}.

We are now ready to prove the main theorem.
\begin{theorem}
Let $M$ be an orientable Poisson manifold of dimension $n$ with a volume form $\mu$ and let $\phi $ be the modular vector field of $M$ with respect to $\mu$. Then for any Poisson geometric module $(W,\{-,-\}_{W})$ and the twisted Poisson geometric module $W_{-\phi }$ of $W$ by $-\phi $, one has the following commutative diagram:\label{thm-poincare-dual-chain}
$$\begin{tikzcd}
\mathfrak{X}^{k}(W) \arrow[d, "\blacktriangle_{W}^{k}"'] \arrow[rr, "\delta_{W}^{k}"] &  & \mathfrak{X}^{k+1}(W) \arrow[d, "\blacktriangle_{W}^{k+1}"] \\
W\otimes_{C}\Omega^{n-k}(M) \arrow[rr, "\partial_{n-k}^{W_{-\phi}}"]       &  & W\otimes_{C}\Omega^{n-k-1}(M)                  
\end{tikzcd}$$
where $\blacktriangle_{W}^{k}=(-1)^{k(k+1)/2}\star_{W}^{k}$, and $\star_{W}^{k}$ is the isomorphism in Proposition \ref{prop-starW-iso}.
\end{theorem}
\begin{proof}
It suffices to verify that $\partial_{n-k}^{W_{-\phi}}\star_{W}^{k}=(-1)^{k-1}\star_{W}^{k+1}\delta_{W}^{k}$. From Lemma \ref{lem-twisted-diff}, we see that for any $X\in \mathfrak{X}^{k}(W)$, 
$$\partial_{n-k}^{W_{-\phi}}\star_{W}^{k}(X)=\partial_{n-k}^{W_{-\phi}}\iota_{X}\mu=\partial_{n-k}^{W}\iota_{X}\mu-(\text{id}_{W}\otimes \iota_{\phi})\iota_{X}\mu.$$
Hence $\partial_{n-k}^{W_{-\phi}}\star_{W}^{k}(X)=(-1)^{k}(\iota_{X}\partial-\iota_{\delta_{W}X})\mu-(\text{id}_{W}\otimes \iota_{\phi})\iota_{X}\mu$, by Lemma \ref{lem-inter-chaindiff}. It follows from (\ref{modular-diff-voleq}) that $\iota_{X}\partial(\mu)=\iota_{X}\iota_{\phi}\mu$. Hence one has
$$\partial_{n-k}^{W_{-\phi}}\star_{W}^{k}(X)=(-1)^{k}\iota_{X}\iota_{\phi}\mu+(-1)^{k-1}\iota_{\delta_{W}X}(\mu)-(\text{id}_{W}\otimes \iota_{\phi})\iota_{X}\mu.$$
By applying Proposition \ref{graded-commutat-iota-prop}, $(-1)^{k}\iota_{X}\iota_{\phi}\mu=(\text{id}_{W}\otimes \iota_{\phi})\iota_{X}\mu$, which completes the proof.
\end{proof}
\begin{remark}
This theorem can be seen as the geometric analogue of \cite[Theorem 3.12]{MR4727073}, but it cannot be directly derived from \cite[Theorem 3.12]{MR4727073} because the commutative algebra $C=C^{\infty}(M)$ here is  not a finitely generated smooth algebra, in general. Usually, $\Omega_{\mathbb{R}}^{1}(C^{\infty}(M))$ is not a finitely generated projective $C^{\infty}(M)$-module, and $C^{\infty}(M)$ is not a Noetherian ring.
\end{remark}
\begin{theorem}\label{thm-poincare-dual-homo}
Let $M$ be an orientable Poisson manifold of dimension $n$ with a volume form $\mu$, and let $\phi $ be the modular vector field of $M$ with respect to $\mu$. Then for any Poisson geometric module $W$ and its twisted Poisson module $W_{-\phi}$ by $-\phi$,
$$\text{HP}_{n-k}(M,W_{-\phi})\cong \text{HP}^{k}(M,W ),\forall 0\leq k\leq n.$$
In particular, $\text{HP}_{n-k}(M,W)\cong \text{HP}^{k}(M,W_{\phi}),\forall 0\leq k\leq n$.
\end{theorem}
For the convenience of readers, we summarize Theorem \ref{thm-poincare-dual-chain} and Theorem \ref{thm-poincare-dual-homo} in the special case where the coefficient Poisson geometric module arises from a flat $T^{*}M$-connection.
\begin{theorem}\label{thm-poincare-flatconnver}
Let $M$ be an orientable Poisson manifold of dimension $n$ with a volume form $\mu$ and let $\phi $ be the modular vector field of $M$ with respect to $\mu$. Then for any smooth vector bundle $E$ over $M$ with a flat $T^{*}M$-connection $\nabla$ and its twisted flat $T^{*}M$-connection $\nabla^{-\phi}$ by $-\phi$, one has the following commutative diagram:
$$\begin{tikzcd}
\mathfrak{X}^{k}(\Gamma(E)) \arrow[d, "\blacktriangle_{E}^{k}"'] \arrow[rr, "\delta_{\nabla}^{k}"] &  & \mathfrak{X}^{k+1}(\Gamma(E)) \arrow[d, "\blacktriangle_{E}^{k+1}"] \\
\Gamma(E\otimes \wedge^{n-k}T^{*}M) \arrow[rr, "\partial_{n-k}^{\nabla^{-\phi}}"]       &  & \Gamma(E\otimes \wedge^{n-k-1}T^{*}M)                  
\end{tikzcd}$$
where $\blacktriangle_{E}^{k}=(-1)^{k(k+1)/2}\star_{\Gamma(E)}^{k}$, and $\star_{\Gamma(E)}^{k}$ is the isomorphism in Proposition \ref{prop-starW-iso}. In particular,
$$\text{HP}_{n-k}(M,\nabla^{-\phi})\cong \text{HP}^{k}(M,\nabla ),\forall 0\leq k\leq n;$$
$$\text{HP}_{n-k}(M,\nabla)\cong \text{HP}^{k}(M,\nabla^{\phi} ),\forall 0\leq k\leq n.$$
\end{theorem}
If $M$ is an orientable unimodular Poisson manifold, it is well known that $HP_{n-k}(M)\cong HP^{k}(M),\forall 0\leq k\leq n$ \cite{MR1675117}. This duality can be easily extended as follows. Since $M$ is unimodular, one can choose a volume form preserved by all Hamiltonian flows, say $\lambda$, then the modular vector field with respect to $\lambda$ is zero. It follows from Theorem \ref{thm-poincare-flatconnver} immediately that for any vector bundle $E$ over $M$ with a flat $T^{*}M$-connection $\nabla$, the diagram
$$\begin{tikzcd}
\mathfrak{X}^{k}(\Gamma(E)) \arrow[d, "\blacktriangle_{E}^{k}"'] \arrow[rr, "\delta_{\nabla}^{k}"] &  & \mathfrak{X}^{k+1}(\Gamma(E)) \arrow[d, "\blacktriangle_{E}^{k+1}"] \\
\Gamma(E\otimes \wedge^{n-k}T^{*}M) \arrow[rr, "\partial_{n-k}^{\nabla }"]       &  & \Gamma(E\otimes \wedge^{n-k-1}T^{*}M)               
\end{tikzcd}$$
is commutative. Hence $HP_{n-k}(M,\nabla)\cong HP^{k}(M,\nabla),\forall 0\leq k\leq n$.
\par Now we relate our twisted duality theorem above with the one in \cite{MR1726784}.
\begin{corollary}
[\text{\cite[Corollary 4.6]{MR1726784}}]Let $M$ be an orientable Poisson manifold of dimension $n$ and let $\nabla^{ELW}$ be the flat $T^{*}M$-connection on the canonical line bundle $\wedge^{n}T^{*}M$ defined by $\nabla_{\omega}^{ELW}(s)=\omega\wedge d\iota_{\pi}(s),\forall \omega\in \Omega^{1}(M),s\in \Omega^{n}(M)$. Then 
$$\text{HP}^{k}(M,\nabla^{ELW})\cong \text{HP}_{n-k}(M),\forall 0\leq k\leq n.$$
\end{corollary}
\begin{proof}
Fix a volume form on $M$, say $\mu$. From Proposition \ref{ELW-twistedill-prop}, we see that $\nabla^{ELW}$ is the twisted flat $T^{*}M$-connection of the obvious flat $T^{*}M$-connection $\nabla_{0}$ on $\wedge^{n}T^{*}M$ (defined by $\mu$) by the modular vector field $\phi_{\mu}$ with respect to $\mu$. Now applying Theorem \ref{thm-poincare-flatconnver}, one has $\text{HP}^{k}(M,\nabla^{ELW})=\text{HP}^{k}(M,\nabla_{0}^{\phi_{\mu}})\cong \text{HP}_{n-k}(M,\nabla)\cong \text{HP}_{n-k}(M),\forall 0\leq k\leq n.$
\end{proof}

\end{document}